\documentclass[12pt]{article}
\usepackage{amsmath,amssymb,amsfonts,epsfig,xspace}
\newcommand{\aref}[1]{Appendix~\ref{sec:#1}}
\newcommand{\sref}[1]{\S~\ref{sec:#1}}
\newcommand{\tref}[1]{Theorem~\ref{thm:#1}}
\newcommand{\lref}[1]{Lemma~\ref{lem:#1}}
\newcommand{\pref}[1]{Proposition~\ref{prp:#1}}
\newcommand{\cref}[1]{Corollary~\ref{cor:#1}}
\newcommand{\fref}[1]{Figure~\ref{fig:#1}}

\newcommand{\old}[1]{}

\newcommand{\Li}{\operatorname{Li}}
\renewcommand{\L}{K} 
\newcommand{\coef}{\gamma}

\newcommand{\eps}{\epsilon}
\newcommand{\veps}{\varepsilon}
\newcommand{\N}{\mathbb N}
\newcommand{\Z}{\mathbb Z}
\newcommand{\R}{\mathbb R}
\newcommand{\C}{\mathbb C}
\newcommand{\ddcl}{\overbrace{c\frac{\partial}{\partial c} \cdots c\frac{\partial}{\partial c}}^\ell}
\newcommand{\enc}{\langle N_c\rangle}
\newcommand{\vnc}{\sigma^2(N_c)}
\newcommand{\st}{\sigma\tau}
\newcommand{\sto}{\sigma_1\tau_1}
\newcommand{\stt}{\sigma_2\tau_2}
\newcommand{\stoc}{\sigma_1,\tau_1}
\newcommand{\sttc}{\sigma_2,\tau_2}
\newcommand{\Zst}{Z_{\st}}
\newcommand{\tqsp}{\tau q + \sigma p}
\newcommand{\tqspo}{\tau_1 q + \sigma_1 p}
\newcommand{\tqspt}{\tau_2 q + \sigma_2 p}
\newcommand{\ii}{\imath}
\newtheorem{thm}{Theorem}
\newtheorem{prop}{Proposition}
\newtheorem{lemma}{Lemma}

\newcommand{\qed}{\hspace*{\fill} \ensuremath{\square}}
\newenvironment{proof}{\noindent{\bf Proof:}\hspace*{1em}}{\qed}
\newenvironment{proofm}{\noindent{\bf Proof:}\hspace*{1em}}{}

\begin{document}
\title{\vspace*{-20pt}
\footnotetext{The research leading to this article was conducted in part while the first author was visiting Microsoft.}
\centerline{\hbox{Critical resonance in the non-intersecting lattice path model}}}
\author{
\begin{tabular}{c}
Richard W. Kenyon \\
 \small CNRS UMR 8628\\
 \small Laboratoire de Math\'ematiques\\
 \small Universit\'e Paris-Sud\\
 \small France\\
 \small\texttt{richard.kenyon\char64math.u-psud.fr}
\end{tabular}
\and
\begin{tabular}{c}
David B. Wilson \\
 \small Microsoft Research\\
 \small One Microsoft Way\\
 \small Redmond, Washington\\
 \small U.S.A.\\
 \small\texttt{dbwilson\char64microsoft.com}
\end{tabular}
}
\date{}
\maketitle

\begin{abstract}
We study the phase transition in the honeycomb dimer model (equivalently,
monotone non-intersecting lattice path model).
At the critical point the system has a strong long-range dependence;
in particular, periodic boundary conditions give rise to a
``resonance'' phenomenon, where the partition function and other
properties of the system depend sensitively on the shape of the domain.
\end{abstract}

\section{Introduction}

We study a model of monotone non-intersecting lattice paths in
$\Z^2$.
Applications of this model include random surfaces
\cite{blote-hilhorst:roughening, Spohn, kenyon:dimers}, magnetic flux lines in
superconductors \cite{WH, BCK}, and a number of other physical phenomena
\cite{Nagle, Fisher, dN}, including spin-domain boundaries
of the three-dimensional Ising
model at zero temperature \cite{cerf-kenyon:wulff}.
(See also \cite{popkov-kim-huang-wu:3d,huang-wu-kunz:5v}.)

Let $R_{m,n}$ be a domain consisting of the $m\times n$ rectangle
in $\Z^2$ with periodic boundary conditions ($R_{m,n}$ is a graph on a torus).
On $R_{m,n}$ consider configurations consisting of collections of
vertex-disjoint, monotone
northeast-going lattice paths (or rather loops, since the paths
are required to eventually close up, possibly after winding several times around
the torus: there are no ``free ends'').  See \fref{bij}.
We do not restrict the number of disjoint loops
but rather
give a configuration an energy $E_b N_b+E_c N_c$ where $N_b$ is the total
number
of ``east'' steps of the paths and $N_c$ is the total number of ``north''
steps of
the paths.  The Boltzmann measure $\mu$ at temperature $T$
is the probability measure assigning a configuration a probability
proportional to $e^{-(E_b N_b+E_c N_c)/T}$.  We study these Boltzmann
measures near the critical temperature $T$, which is the temperature at which
$e^{-E_b/T}+e^{-E_c/T}=1$ \cite{kasteleyn:pfaffian}.
Letting $b= e^{-E_b/T}$ and
$c=e^{-E_c/T}$, a configuration
has a probability proportional to $b^{N_b}c^{N_c}$.

This process is equivalent to another well-known model,
the model of dimers (weighted perfect matchings) on the
honeycomb lattice $H$ \cite{Wu68,WH:93}.
Weight the edges of the honeycomb lattice
$a=1$, $b$, or $c$ according to their direction as in \fref{bij}.
See \fref{bij} for an illustration of the weight-preserving bijection
between dimer configurations and lattice paths.

\begin{figure}[htbp]
\centerline{\epsfig{figure=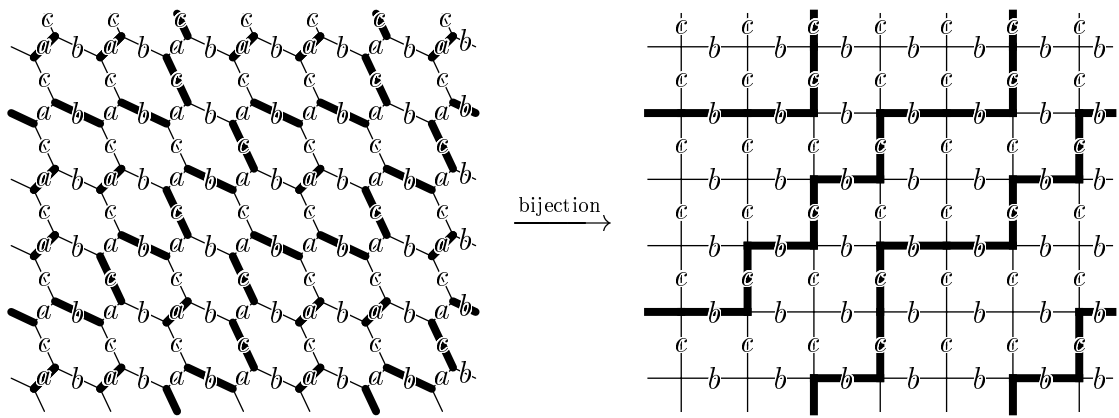}}
\caption{
The classical bijection \cite{Wu68} between perfect matchings of the hexagonal
lattice and non-intersecting north-east lattice paths, shown for a
$m\times n = 7\times 6$ region.  Edges of type `$a$' in the hexagonal
lattice $H_{m,n}$ are contracted to obtain $R_{m,n}$.  The dimer
configuration has weight $a^{16} b^{14} c^{12}$, as does the lattice
path configuration (when $a=1$ or else one introduces a factor of $a$
for each vertex not in a lattice path).}
\label{fig:bij}
\end{figure}

Let $H_{m,n}$ be a finite graph which is a quotient of $H$ by the
horizontal and vertical translations of length $m,n$ respectively
as in \fref{bij}.
Let $\mu_{m,n}=\mu_{m,n}(a,b,c)$ be the Boltzmann measure on
dimer configurations on the toroidal graph $H_{m,n}$.
This measure assigns a dimer configuration a probability proportional
to $a^{N_a}b^{N_b}c^{N_c}$, where there are $N_a$ edges of weight
$a$, $N_b$ edges of weight $b$, and $N_c$ edges of weight $c$.
The exact probability is $a^{N_a}b^{N_b}c^{N_c}/Z$, where
the normalizing constant $Z$ is called the partition function
of the system.

As $m,n\to\infty$ the measures $\mu_{m,n}(a,b,c)$ have a unique limiting
Gibbs measure $\mu(a,b,c)$ \cite{Sheffield, cohn-kenyon-propp:variational}.
The measure $\mu(a,b,c)$ is a measure on dimer configurations on $H$,
and is well-understood for all $a,b,c$  \cite{cohn-kenyon-propp:variational}.
As $a,b,c$ vary the measure $\mu(a,b,c)$ undergoes a phase transition
--the ``solid-liquid'' transition-- when
one of $a$, $b$, or $c$ is equal to the sum of the other two,
for example $a=b+c$.  When $a\geq b+c$ the system is frozen:
with probability $1$ only edges of weight $a$ are present.  See \fref{melt}.
When $a,b,c$ satisfy the strict triangle inequality (each is less
than the sum of the others) dimers of all types are
present in a typical configuration of $\mu(a,b,c)$.
Likewise in the lattice path model, when $a\geq b+c$ there are no lattice paths.
When $a,b,c$ satisfy the strict triangle inequality,
the lattice paths in a configuration are dense, that is, on average they lie
within a constant distance of one another.

\begin{figure}[htbp]
\centerline{\epsfig{figure=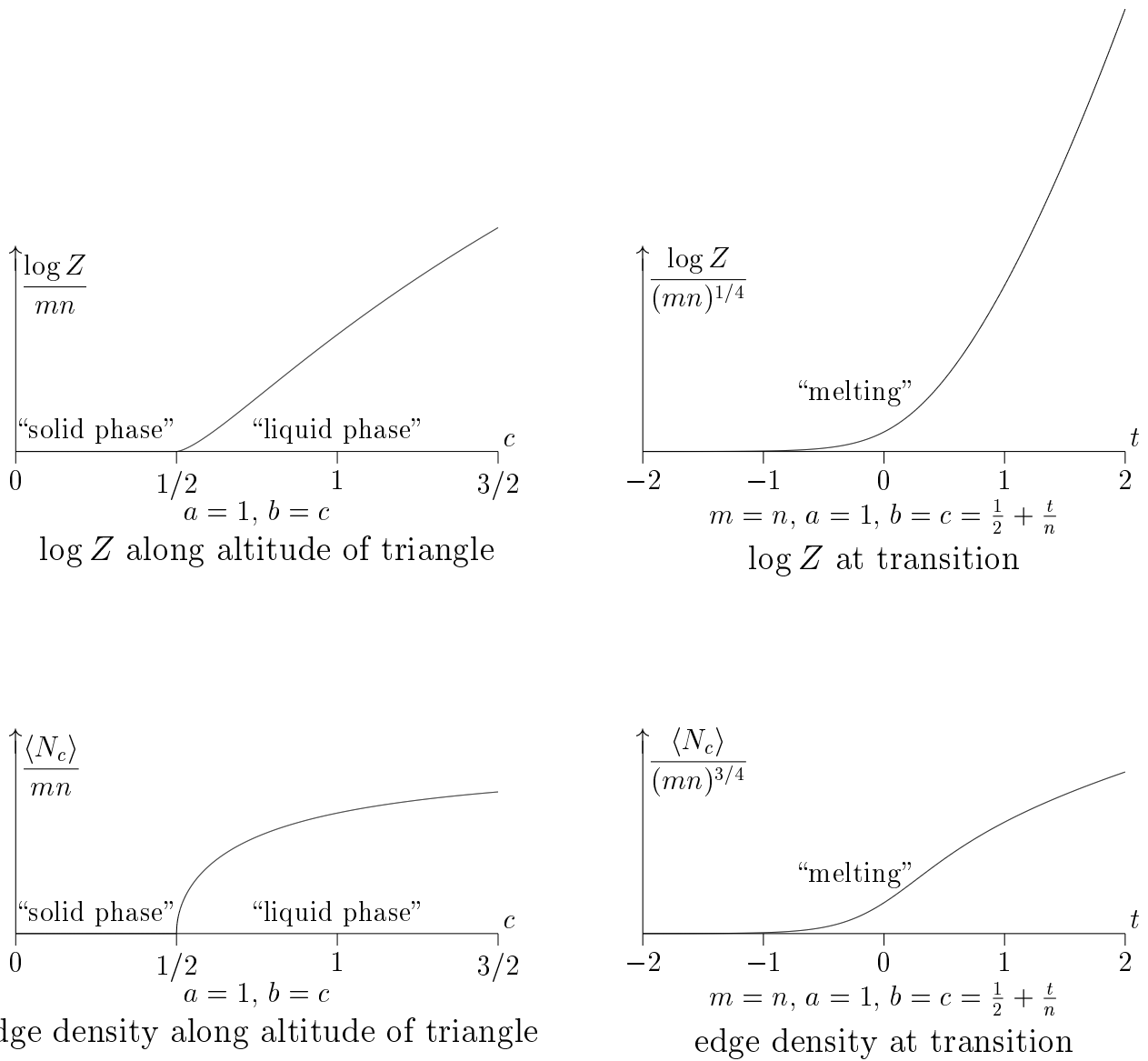}}
\caption{
The partition function is shown in the upper panels, and the expected
number of edges of type `$c$' is shown in the lower panels.  The solid
and liquid phases are shown in the left panels, and the melting
transition between between these phases is shown in the right panels.
The formulas on the left are derived in
\cite{kasteleyn:pfaffian,Wu67,cohn-kenyon-propp:variational} and are not
used here.  The formulas on the right follow from \tref{gaussian}.
The melting transition depends quite sensitively on the aspect ratio
of the region (see \fref{spike} and \tref{intro}).}
\label{fig:melt}
\end{figure}

The finite-volume measures $\mu_{m,n}(a,b,c)$ are
less well understood near $a=b+c$.  When $a=b+c$ the lattice
paths exist but are spread out;
as we shall see,
the average distance between strands is on the order of
the square root of the system size.
The ratio $n/m$ imposes fairly rigid entropic constraints
on the way these loops can join up.  Surprisingly these constraints
become stronger with increasing system size, so that the scaling
limit has nontrivial structure.

Our main result is a computation of the
partition function as a function of $a,b,c,m,n$, for parameters
near $a=b+c$.  
The partition function $Z$ (normalized by $(\text{area})^{1/4}$)
is largest when $nb/(mc)$ is rational with small numerator and denominator.

For example we have
\begin{thm}\label{thm:intro}
When $a=1,b=c=1/2$ and $n/m=p/q$ in lowest terms we have
$$\log Z = \frac{(mn)^{1/4}\zeta(3/2)(1-2^{-1/2})}{2\sqrt{\pi}(pq)^{3/4}}(1+o(1)).$$
as $m$ and $n$ tend to $\infty$ while $p$ and $q$ remain fixed
(here $\zeta$ is the Riemann zeta function).  Furthermore
the number $N_c$ of `$c$'-type edges tends to a Gaussian with expectation
$$\enc=\frac{(mn)^{3/4}\zeta(1/2)(1-2^{1/2})}{2\sqrt{\pi}(pq)^{1/4}}(1+o(1))$$
and variance
$$\vnc=\frac{(mn)^{5/4}(pq)^{1/4}\zeta(-1/2)(1-2^{3/2})}{2\sqrt{\pi}}(1+o(1)).$$
\end{thm}
(Both $\zeta(1/2)$ and $\zeta(-1/2)$ are negative.)

We see that the system greatly
prefers domains of {\it rational} modulus over those with irrational
modulus.  Here by rational we mean, a lattice path with
average slope $c/b$ will close up after winding a small number of
times around the torus.  In such a case the number of loops can be large,
whereas in an irrational case each loop must wind many times
around before closing up (unless it pays a large entropic cost).

One can think of the partition function $Z$,
taken as a function of $m$ and $n$, as an indicator of
rationality of $n/m$.
See \fref{spike} which plots $\log Z / (\text{area})^{1/4}$
and $\enc / (\text{area})^{3/4}$ as a function
of $\log(n/m)$ near $\text{area}=10^7$ for $a=1$ and $b=c=\frac12$.

\begin{figure}[htbp]
\centerline{\epsfig{figure=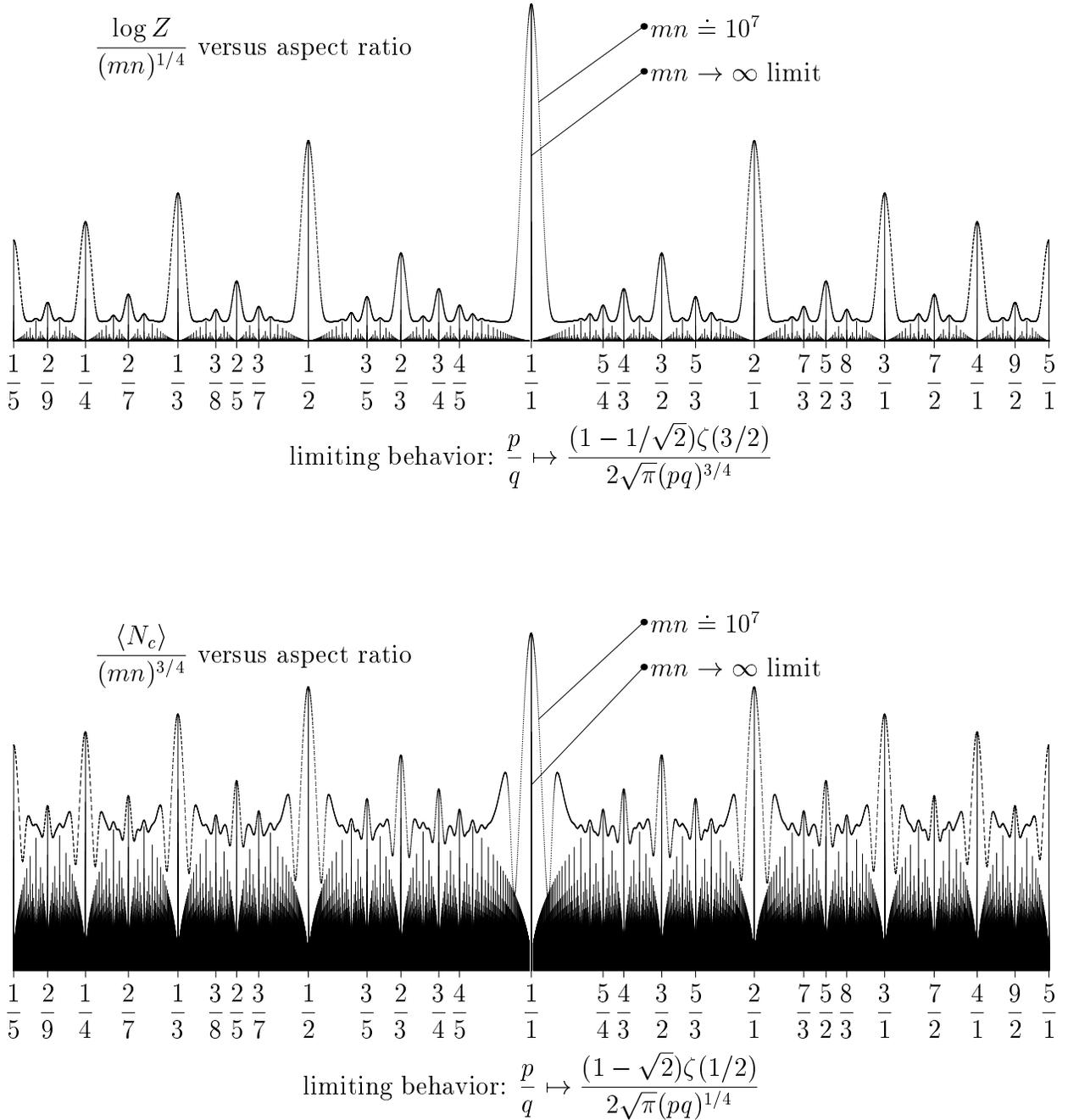}}
\caption{Resonant spikes in $\log Z$ (upper panel) and in the edge density (lower panel) for very large regions when $a=1$ and $b=c=1/2$.  The aspect ratio is plotted on a log scale.}
\label{fig:spike}
\end{figure}

This is not the only interesting phenomenon in this model.
For a rational domain, there is
a non-trivial behavior as we vary $a,b,c$ away from the critical
point.  Letting $A=(c/(a-b))^n$, the partition function as a function of $A$
has an infinite number of nonanalyticities (in the large-$n$ limit)
which correspond to abrupt changes in the winding number of
curves in a typical configuration.  That is, the curves ``ratchet'' at
well-defined values as $A$ increases:
see \fref{Zsig} and \sref{crossover}. We did not prove (although we believe)
that at a typical value of $A$ the curves are in a well-defined
integer homology class 
(i.e., have well-defined winding number around the torus),
and this homology class changes at discrete values of $A$.
We prove only that the $\Z_2$ homology class (winding number modulo $2$) 
changes at these well-defined points.
We also show that, away from these transition points, 
the number of edges is a Gaussian;
at the transition it is a mixture of two Gaussians, coming from a Gaussian
for each homology class.

It would be very interesting to study this same model on higher-genus surfaces.
On  higher-genus surfaces (with translation structures having
conical singularities) it would be very useful to
be able to detect ``rationality'', in the form of bands of parallel
closed geodesics: this is an important problem in
billiards \cite{Tabachnikov}.  Moreover the ratcheting phenomenon
must be significantly richer in the presence of a non-abelian
fundamental group.

\section{Review of $\log Z$ and the distribution of $N_c$} \label{sec:logZ}

The partition function $Z$ is
$$Z=\sum_{\text{configurations}}a^{N_a}b^{N_b}c^{N_c}.$$
It contains all the information about the distribution of the
total number of edges of each type.
We can view $Z$ as a polynomial in $c$, where the
coefficient of $c^{N_c}$ is the weighted sum of configurations which
contain $N_c$ edges of type `$c$'.  Thus the expected number of edges of
type `$c$' is
$\enc=\frac{c}{Z}\frac{\partial Z}{\partial c}$.  Similarly,
the $\ell$th moment of the number of edges of type `$c$' is given by
\newcommand{\M}[1]{\langle N_c^{#1}\rangle}
$$\langle N_c^\ell\rangle=\frac{1}{Z} \ddcl Z.$$
The $\ell$th {\it cumulant} $\L_\ell$ is defined by
$$\L_\ell = \ddcl \log Z,$$
which we will in effect estimate later.
Here we recall several properties of cumulants \cite{Comtet}.
Of course $\M{1}=\L_1$.
The variance in the number of edges of type `$c$' is $\L_2$:
\begin{align*}
c\frac{\partial}{\partial c}c\frac{\partial\log Z}{\partial c}
 &= c\frac{\partial}{\partial c}\frac{c}{Z}\frac{\partial Z}{\partial c} \\
 &= \frac{c}{Z}\frac{\partial}{\partial c}c\frac{\partial Z}{\partial c}
 - \frac{c}{Z^2} \frac{\partial Z}{\partial c} c\frac{\partial Z}{\partial c} \\
 &= \langle N_c^2\rangle - \enc^2 = \vnc.
\end{align*}

Later we will use higher moments to show that the number of `$c$'-type edges tends to a Gaussian, and to this end we express these moments in terms of the $\L_\ell$'s.
Since $$\M{\ell+1} = \frac{1}{Z} \frac{c \partial}{\partial c} (Z \M{\ell}),$$
when $\M{\ell}$ is expressed as a polynomial
in the variables $\L_1,\ldots,\L_\ell$,
we may calculate $\M{\ell+1}$ from $\M{\ell}$ by replacing each monomial
$\L_{i_1} \L_{i_2} \cdots \L_{i_k}$ of $\M{\ell}$ with
$$\L_1 \L_{i_1} \L_{i_2} \cdots \L_{i_k} + \L_{i_1+1} \L_{i_2}
\cdots \L_{i_k} + \L_{i_1} \L_{i_2+1} \cdots \L_{i_k} + \cdots + \L_{i_1} \L_{i_2}
\cdots \L_{i_k+1}.$$  Thus for example we have
\begin{align}
\label{Bellpoly}
\begin{split}
\M{1} =& \L_1\\
\M{2} =& \L_2+\L_1^2\\
\M{3} =& \L_3 + 3 \L_2 \L_1 + \L_1^3\\
\M{4} =& \L_4 + 4 \L_3 \L_1 + 3 \L_2^2 + 6 \L_2 \L_1^2 + \L_1^4\\
\M{5} =& \L_5 + 5 \L_4 \L_1 + 10 \L_3 \L_2 + 10 \L_3 \L_1^2 + 15 \L_2^2 \L_1 + 10 \L_2 \L_1^3 + \L_1^5.
\end{split}
\end{align}
These polynomials $\M{j}=Y_j(\L_1,\L_2,\dots)$ are the
{\it complete Bell polynomials} \cite{Comtet}.
We see that $\M{\ell}$ contains a monomial for each partition of $\ell$,
and the coefficient associated with partition with distinct part sizes $s_1 > s_2 >
\cdots > s_k$ and $r_i$ parts of size $s_i$ is
\begin{equation}
\frac{\ell!}{s_1!^{r_1} r_1! \cdots s_k!^{r_k} r_k!}. \label{partcoef}
\end{equation}

It will be more useful to work with moments about the mean rather
than moments about the origin.  Note that if we replace $Z$ with $Z^*
= Z c^{-\mu}$, then the above derivation shows us how to express
$\langle (N_c-\mu)^\ell\rangle$ in terms of the $\L^*_\ell$'s defined by
$$\L^*_\ell = \ddcl \log (Z c^{-\mu})
= \L_\ell + \ddcl \log (c^{-\mu}).$$
As $\L^*_1 = \L_1 -\mu$ and $\L^*_\ell = \L_\ell$ for $\ell>1$, upon
substituting $\mu = \enc = \L_1$ we see that the $\ell$th
moment of $N_c$ about the mean may be obtained from the above
expressions for $\M{\ell}$ by deleting all monomials that contain the
variable $\L_1$.

\section{Product form of the partition function}

We compute an expression for the partition
function as a function of $a,b,$ and $c$.
We are interested in approximating $Z$ to within
$1+o(1)$ multiplicative errors when $a=1$, $b,c\in(0,1)$ and $b+c$ 
is close to $a$.
The interesting range is when $(c/(a-b))^n$ is of constant order, that is,
$a-b-c=O(1/n)$.  We define $A=(c/(a-b))^n$ as the natural parameter
measuring proximity to the critical point.

In what follows we always set $a=1$, although we keep using $a$
for notational convenience.

Recall $H_{m,n}$, the $m\times n$ hexagonal toroidal graph shown in
\fref{bij}.  By Kasteleyn \cite{kasteleyn:pfaffian} (see also
\cite{tesler:pfaffian,galluccio-loebl:pfaffian,regge-zecchina:pfaffian,lu-wu:dimer1,lu-wu:dimer2}
for extensions and further developments), the partition function
$Z=Z(a,b,c)$ for dimer coverings of $H_{m,n}$
is a sum of four expressions,
\begin{equation}\label{Z4}
Z = \frac12(-Z_{00}+Z_{01}+Z_{10}+Z_{11}),
\end{equation}
where $\Zst$ is the determinant of a signed version of the adjacency
matrix of $H_{m,n}$,
and counts dimer coverings with a sign according to the homology class
(in $H_1(\text{torus},\Z_2)\cong \Z_2^2$)
of the corresponding system of loops, as follows.
\newcommand{\vepsx}{\veps_{\hat x}}
\newcommand{\vepsy}{\veps_{\hat y}}
Let $N(\vepsx,\vepsy)$ denote the total weight of dimer coverings whose corresponding
loops have $\vepsx\bmod 2$ crossings of the line $x=0$
and $\vepsy\bmod 2$ crossings of the line $y=0$.
Each $\Zst$ is a linear combination
of the $N(\vepsx,\vepsy)$ with coefficients
$\pm1$ as follows:
\begin{equation}\label{Ztable}\begin{array}{rcccc}
&N(0,0)&N(1,0)&N(0,1)&N(1,1)\\
Z_{00}&+1&-1&-1&-1\\
Z_{10}&+1&+1&-1&+1\\
Z_{01}&+1&-1&+1&+1\\
Z_{11}&+1&+1&+1&-1\makebox[0pt][l]{.}\\
\end{array}
\end{equation}
Note three important facts, which follow from this table:
\begin{prop}\label{prp:2Zsums}
The sum of any two of $-Z_{00},Z_{01},Z_{10},Z_{11}$
has only nonnegative coefficients.
The difference between any two
of $-Z_{00},Z_{01},Z_{10},Z_{11}$ is bounded by the sum of the other two.
The difference between the coefficients of $a^{\alpha} b^{\beta} c^{\gamma}$
in any two of $-Z_{00},Z_{01},Z_{10},Z_{11}$ is bounded by the sum of the
coefficients of $a^{\alpha} b^{\beta} c^{\gamma}$ in the other two.
\end{prop}

Kasteleyn \cite{kasteleyn:pfaffian} evaluated the determinants $\Zst$ by
multiplying eigenvalues obtained through Fourier analysis, giving
\enlargethispage*{\baselineskip}
\newcommand{\ms}{(-1)^\sigma}
\newcommand{\mt}{(-1)^\tau}
\newcommand{\mst}{\veps_{\st}}
\newcommand{\msto}{\veps_{\sto}}
\newcommand{\mstt}{\veps_{\stt}}
\begin{align}\nonumber
\Zst
\nonumber
  &= \prod_{\substack{(-z)^m=\ms\\(-w)^n=\mt}} [ a+b z+c w ] \\
\nonumber
  &= \prod_{(-z)^m=\ms} [ (a+bz)^n -\mt c^n  ]\\
\nonumber
  &= \prod_{(-z)^m=\ms}(a+bz)^n
     \prod_{(-z)^m=\ms}\left[1 -\mt \left(\frac{c}{a+bz}\right)^n\right] \\
\label{exact}  &= (a^m -\ms b^m)^n
     \prod_{(-z)^m=\ms}\left[1 - \mt \left(\frac{c}{a+bz}\right)^n\right].
\end{align}

We can ignore the $b^m$ term which is exponentially smaller than $a^m$.
When $b+c$ is close to $a$, unless $z$ is close to $-1$,
$|a+bz|$ will be greater than $c$;
in particular $|a+bz|^n$
is exponentially larger than $c^n$.
So we can ignore the factors in the products for which
$z$ is not close to $-1$.  We will expand the remaining factors near $z=-1$.
As $z$ is a root of unity, let $z=z_k=-e^{i\theta_k}$, where
$\theta_k = 2\pi k/m$ for $k \in \Z_m+\sigma/2$.  Of course $k\equiv k+m$,
so when doing series expansions we can take $-m/2<k\leq m/2$.

Define $r_k,\phi_k$ by
\begin{align*}
1\pm \left(\frac{c}{a+b z_k}\right)^n &=
1\pm \left(\frac{c}{a-b}\right)^n r_k e^{i\phi_k},
\end{align*}
so that
 $$r_k=(a-b)^{n}|a+bz_k|^{-n} \ \ \ \ \text{and}\ \ \ \ 
  \phi_k=\arg(a+bz_k)^{-n}.$$
\enlargethispage*{\baselineskip}
We make the simplifying assumptions $b/a=\Theta(1)$, $1-b/a=\Theta(1)$, and $n=\Theta(m)$.  Then
$$r_k=\exp\left(-n\theta_k^2\frac{a b}{2(a-b)^2}+O(n\theta_k^4)\right)=
\exp(-\eps k^2+O(k^4/m^3))$$
where we have defined $$\eps = \frac{2\pi^2 n a b}{m^2(a-b)^2} = O(1/m),$$
and
$$\phi_k = \frac{n\theta_k b}{a-b} +n O(\theta_k^3)
         = \phi k + O(k^3/m^2)$$
where we define $$\phi = \frac{2\pi n b}{m(a-b)}.$$
Letting $A=\left(\frac{c}{a-b}\right)^n$, we have
\begin{equation}\label{pf}
 \Zst = (a^m -\ms b^m)^n \prod_{k\in\Z_m+\sigma/2} (1 -\mt A r_k e^{i \phi_k}).
\end{equation}
In logarithmic form we can write
\begin{equation}\label{sumform}
\log \Zst = - \sum_{k\in\Z_m+\sigma/2}
   \Li_1\big((-1)^\tau A r_k e^{i \phi_k}\big)
 + \underbrace{n\log(a^m-\ms b^m)}_{\text{negligible}}
\end{equation}
where the polylogarithm function $\Li_\nu$ is defined by
$\Li_\nu(z)=\sum_{n=1}^{\infty}\frac{z^n}{n^\nu}$ for $|z|<1$ and by
analytic continuation elsewhere (see
\aref{polylogs} for background on polylogarithms).
Here the second term is essentially zero,
and the terms in the summation become negligible when
$|k|$ is larger than $\Theta(\sqrt{m})$.

In this expression for $\log \Zst$ we need only keep track of
its real part: from \eqref{pf} we see that (since $b<a$)
$\Zst$ is real and nonnegative if $A\leq1$,
and real and nonnegative if $A>1$ except for $Z_{00}$ which is strictly
negative.  In particular when $A>1$ the expression \eqref{Z4}
is a sum of nonnegative terms.
When $A\leq 1$ we shall see that $Z_{00}$
is negligible compared to $Z$ so its sign is irrelevant.

More generally we find, using \eqref{sumform},
$z\frac{d}{dz}\Li_\nu(z)=\Li_{\nu-1}(z)$, and
$(c\partial/\partial c)A=nA$, that
\begin{equation}\label{sumform2}
\begin{split}
\ddcl \log(\Zst) =
- n^{\ell}\sum_{k\in\Z_m+\sigma/2} &\Li_{1-\ell}\big(\mt A r_k e^{i \phi_k}\big)
\\&\text{($+\underbrace{n\log(a^m-\ms b^m)}_{\text{negligible}}$ if $\ell=0$).}
\end{split}
\end{equation}

\section{Rational tori} \label{sec:rational}

The expressions \eqref{sumform2}
are non-trivial to evaluate,
mostly because they are describing behavior which depends sensitively
on the parameters defining the system.
In this section we compute the asymptotics of \eqref{sumform2}
for ``nearly rational'' domains.

We consider the toroidal hexagonal graph $H_{m,n}$
to be ``nearly rational'' when
$\phi/(2\pi)=\frac{nb}{m(a-b)}$ is
close to a simple rational $p/q$, where $p$ and $q$ are relatively
prime integers.
(Note that this depends not only on $m,n$ but also on $a,b,c$.)
We keep $p$ and $q$ fixed as
$m$ and $n$ tend to infinity, and by ``close'' we mean that
$nb/(mc)-p/q$ is not large compared to $1/\sqrt{qn}$.
We introduce the parameter $\alpha$ to measure the closeness
of $nb/(mc)$ to $p/q$:
we define $\alpha$ by
$$\frac{\phi}{2\pi} = \frac{p}{q} (1+\alpha W),$$ 
where $W= \sqrt{q \eps}/(\pi p)=\Theta(\frac1{\sqrt{qn}})$.
The interesting range is when $\alpha$ is of constant order.

We will determine the asymptotic shape of the resonant peaks in \fref{spike}
as functions of $\alpha$ and $A$.

\subsection{Spokes, spirals, and clouds}\label{sec:spirals}

For fixed $p$ and $q$, as the area $mn$ gets large, the terms $1\pm A r_k
e^{i \phi_k}\approx 1\pm Ar_ke^{2\pi i(p/q)k}$ accumulate on $q$ different spokes (or radii) of a
circle with radius $A$ centered at $1$.  If $\phi/(2\pi)$
is only approximately $p/q$, then each of these spokes becomes
a spiral, which spirals out from $1$ when $k$ is negative and
increasing, and then spirals back in towards $1$ when $k$ is positive
and increasing.  If $\phi/(2\pi)$ is far from a simple rational,
then the terms $1\pm A r_k e^{i \phi_k}$ form a cloud within the disk
of radius $A$ centered at $1$.  In our analysis for nearly
rational tori, we will assume that $\phi/(2\pi)$ is sufficiently
close to a simple enough rational $p/q$ that the terms $1\pm A r_k
e^{i \phi_k}$ form what appear to be $q$ continuous spokes or spirals
(in a sense we define more precisely below).  See \fref{spirals}.

It is useful to re-express \eqref{sumform2} to reflect the presence of
the $q$ spirals.  Since $\phi\approx 2\pi\frac{p}{q}$, every $q$th
term lies on a given spiral, so we break the sum apart into $q$ different
sums, one for each spiral.
  Of course $q$ may not evenly divide $m$.
The most convenient way to re-express the
 summation is
$$ \sum_{k\in\Z_m+\sigma/2} f(k) =
 \frac{1}{q} \sum_{j\in\Z_q+\sigma/2} \,\,\sum_{u\in\Z_m} f(j+qu).$$
\begin{figure}[t]
\centerline{\epsfig{figure=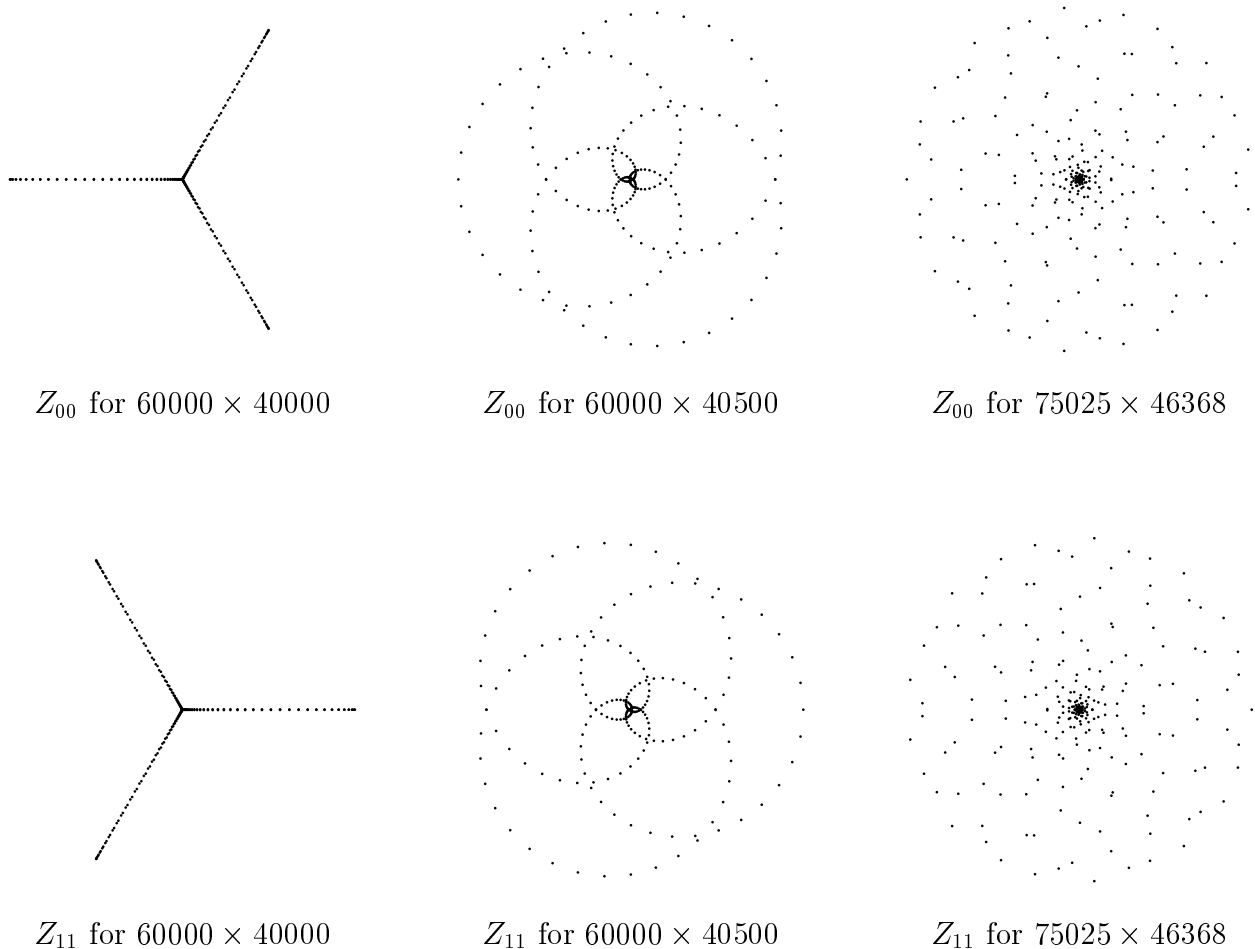}}
\caption{The multiplicands for $Z_{00}$ and $Z_{11}$ for a domain whose aspect ratio is (a) simple rational (b) nearly simple rational (c) far from simple rational.  The multiplicands are complex and accumulate towards $1$.}
\label{fig:spirals}
\end{figure}
Since spirals are continuous objects rather than discrete sets of points,
we wish to approximate $\sum_{u\in\Z_m}$ with $\int_0^m\,du$.  
To this end we subtract $2\pi p u$ from $\phi_{j+qu}$; 
then the angle ($\bmod\ 2\pi$)
is unchanged for integer $u$, while for continuous $u$ it is slowly varying
so we can hope that the integral approximates the sum.
Thus we re-express \eqref{sumform2} as
(ignoring the negligible $O(n b^m)$ error term when $\ell=0$)
\begin{align}
\ddcl \log(\Zst) &\doteq
\begin{aligned}[t]
- n^{\ell} \frac{1}{q} \sum_{j\in\Z_q+\sigma/2} \,\,\sum_{u\in\Z_m}&\Li_{1-\ell}\big(\mt A r_{j+qu} e^{i (\phi_{j+qu}-2\pi p u)}\big)
\end{aligned} \nonumber\\
\intertext{
Later we will quantify the error introduced by approximating
these sums with integrals,
and show it to be insignificant
so long as the points appear to line up on $q$ spirals which miss the
singularity.
For now we proceed with the integral
approximation and simplify it:}
 &\approx
- n^{\ell} \frac{1}{q} \sum_{j\in\Z_q+\sigma/2} \,\,\int_0^m \Li_{1-\ell}\big(\mt A r_{j+qu} e^{i (\phi_{j+qu}-2\pi p u)}\big) \,du \label{intapprox}\\
\intertext{changing variables to $k=j+qu$ and using the fact that the integrand is periodic,}
&=
- n^{\ell} \frac{1}{q^2} \int_0^{qm} \sum_{j\in\Z_q+\sigma/2} \Li_{1-\ell}\big(\mt A r_k e^{i (\phi_k-2\pi k p/q)} e^{2\pi i j p/q)}\big) \,dk \nonumber\\
\intertext{using the replication formula
$ \frac1q \sum_{\omega^q=1} \Li_\nu(\omega z) = \frac{1}{q^\nu}\Li_\nu(z^q)$
and the fact that $\gcd(p,q)=1$,}
&=
- (qn)^{\ell} \frac{1}{q^2} \int_0^{qm} \Li_{1-\ell}\big((-1)^{\tau q} A^q r_k^q e^{i (\phi_k -2\pi k p/q)q} e^{2\pi i (\sigma/2)p}\big) \,dk \nonumber\\
\intertext{and again using the periodicity of $\phi_k$ and $r_k$,}
&=
- \frac{(qn)^{\ell}}{q} \int_{-m/2}^{m/2} \Li_{1-\ell}\big((-1)^{\tau q+\sigma p} A^q r_k^q e^{i (\phi_k -2\pi k p/q)q} \big) \,dk, \label{intform}\\
\intertext{which, as a function of $\sigma$ and $\tau$,
only depends upon the parity of $\tqsp$.
To obtain the next formula we substitute the estimate $r_k e^{i \phi_k}
\approx e^{-\eps (k \bmod m)^2} e^{i \phi (k \bmod m)}$ (where we take
$k \bmod m$ to lie between $-m/2$ and $m/2$).
Doing this substitution introduces an error,
but we postpone the error analysis until later.  Note that we can substitute
this approximation for $r_k e^{i\phi_k}$ either in \eqref{intform}
or just prior to the integral approximation \eqref{intapprox},
since the only property of $r_k e^{i\phi_k}$ that we used in the
intervening steps is that it is periodic in $k$ with period $m$.}
&\approx
- \frac{(qn)^{\ell}}{q} \int_{-m/2}^{m/2} \Li_{1-\ell}\big((-1)^{\tau q+\sigma p} A^q e^{-q\eps k^2} e^{i (\phi -2\pi p/q)q k} \big) \,dk. \label{intform2}\\
\intertext{We extend the range of
integration to all of $\R$, introducing a negligible error,}
&\approx
- \frac{(qn)^{\ell}}{q} \int_{-\infty}^{\infty} \Li_{1-\ell}\big((-1)^{\tau q+\sigma p} A^q e^{-q\eps k^2} e^{i (\phi -2\pi p/q)q k} \big) \,dk. \label{intform3}\\
\intertext{To measure the closeness of $\phi/(2\pi)$ to $p/q$, we
define $\alpha$ so that $\phi/(2\pi) = p/q (1+\alpha W)$, where we define
$W= \sqrt{q \eps}/(\pi p)$.  Then $(\phi - 2\pi p/q) q k = 2\pi p \alpha W k = 2 \alpha\sqrt{q\eps} k$.  We change variables to $x=\sqrt{q\eps} k$ to obtain}
 &= -\frac{(qn)^\ell}{\sqrt{q^3\eps}} \int_{-\infty}^{\infty} \Li_{1-\ell}\big((-1)^{\tau q+\sigma p} A^q e^{2 i \alpha x - x^2}\big) \,dx. \label{ddclZ}
\end{align}

The reader may wonder about the apparent asymmetry in
equation~\eqref{ddclZ} (when $\ell=0$), e.g.\ why 
does $q$ appear but not $p$, while
$Z$ is symmetrical with respect to width and height?  But when
$b+c\approx a\approx 1$ and $(nb)/(mc)\approx p/q$ we have
$$W= \frac{\sqrt{q \eps}}{\pi p} \approx \sqrt{\frac{2}{pmc}} \approx \sqrt{\frac{2}{qnb}} \approx \frac{\sqrt{2}}{(p q m n b c)^{1/4}},$$
\begin{equation}
\frac{1}{\sqrt{q^3\eps}} = \frac{1}{\pi p q W} \approx \frac{(m n b c)^{1/4}}{\pi \sqrt{2} (p q)^{3/4}}, \label{q3eps}
\end{equation}
$$ \log A^q = n q \log\frac{c}{a-b} \approx (b+c-a) \frac{n q}{c} \approx (b+c-a) \frac{m p}{b},$$
$$ (qn)^\ell \approx (pqmnc/b)^{\ell/2}.$$
The asymmetry between $b$ and $c$ when $\ell>0$
should not be unexpected since we differentiated with respect to $c$
rather than $b$.

Referring back to the spirals in \fref{spirals}, $q$ counts the number
of spirals, the parameter $A$ measures the radius of the spirals,
$\sqrt{\eps}$ is a measure of how far apart the points are on the
spiral, and $\alpha$ is a measure of the ``spirality''.  $\alpha$ is the
right parameter against which to plot the shape of the spikes, making
$W$ a measure if their width.  When $\alpha=0$ the spirals are spokes,
when $\alpha$ gets too large (for a given $\sqrt{\eps}$) the spirals
break up into a cloud, by which time the integral approximation
\eqref{intapprox} breaks down.

\subsection{Error analysis}\label{sec:simple}

In the interest of simplicity, we only consider the case when
$p$, $q$, $A$, and $\alpha$ are held fixed while $m\rightarrow\infty$
and $n\rightarrow\infty$.  Most of the interesting behavior already
shows up in this case.  It is also quite interesting to ask how much
these parameters can vary (e.g.\ can $p$ and $q$ be as large as
$n^{1/3}$?), but we do not pursue that in this article.

\begin{lemma} \label{lem:ddclz}
When we fix $p$, $q$, $A$, and $\alpha$ while $m\rightarrow\infty$ and
$n\rightarrow\infty$, we have
$$ \frac{\sqrt{q^3\eps}}{(qn)^\ell}\ddcl \log(\Zst) =
-\int_{-\infty}^{\infty} \Li_{1-\ell}\big((-1)^{\tqsp} A^q 
e^{2 i \alpha x - x^2}\big) \,dx + o(1),$$
provided that the curve $(-1)^{\tqsp} A^q e^{2 i \alpha x - x^2}$ 
does not contain the point $1$.  If the curve does contain 
$1$, then the convergence for $\ell=0$ 
is still valid provided that no multiplicand of $\Zst$ is closer than 
$e^{-o(\sqrt{n})}$ to $0$, in which case the right-hand side is merely an upper bound.
\end{lemma}

These integrals are explicitly evaluated in \sref{Zxint}.

\begin{proof}
Much of the proof has already been given in \sref{spirals}, what we have left 
to do is justify the approximations that we made in \eqref{intapprox} and 
\eqref{intform2}.
For this error analysis we do the
$r_k e^{i\phi_k} \approx e^{-\eps k^2 + i\phi k}$ substitution before the
integral approximation.  With $k=j+qu$ and $j\in\Z_q+\sigma/2$,
\begin{align*}
 \Li_{1-\ell}\big(\mt A r_{j+qu} e^{i (\phi_{j+qu}-2\pi p u)}\big)
&= 
 \Li_{1-\ell}\big(\mt e^{2\pi i j p/q} A e^{i (\phi-2\pi p/q)k - \eps k^2+ O(k^3/m^2)}\big)\\
&= \Li_{1-\ell}\big(\mt e^{2\pi i j p/q} A e^{2 i \alpha\sqrt{\eps/q}\, k - \eps k^2 + O(k^3/m^2)}\big)
\intertext{using $\phi-2\pi p/q = 2\pi(p/q)\alpha W=2\alpha\sqrt{\eps/q}$.
Next we use the fact that for integer $\ell\geq0$,
$\frac{d}{dz}\Li_{1-\ell}(e^z)=\Li_{-\ell}(e^z)$ is a rational function of
$e^z$ with
a pole at $e^z=1$ but which is bounded outside a neighborhood of this pole,}
&= \Li_{1-\ell}\big(\mt e^{2\pi i j p/q} A e^{2 i \alpha\sqrt{\eps/q}\, k - \eps k^2}\big)+O(k^3/m^2)
\end{align*}
which is valid as long as
$\mt e^{2\pi i j p/q} A e^{2 i \alpha\sqrt{\eps/q}\, k - \eps
k^2}$ lies outside a neighborhood of $1$, and the $O(k^3/m^2)$ error
term is much smaller than the radius of this neighborhood.

It is not hard to see that the $q$ curves
$\mt e^{2\pi i j p/q} A e^{2 i \alpha\sqrt{\eps/q}\, k - \eps k^2}$
are bounded away from $1$ if and only if the curve
$(-1)^{\tqsp} A^q e^{2 i \alpha x - x^2}$ avoids $1$.

Adding up the errors $O(k^3/m^2)$ over the range
$|k|<\Theta(\sqrt{n})$ gives $O(1)$.  When $|k|\gg\Theta(\sqrt{n})$,
both $r_k$ and its approximation $e^{-\eps k^2}$ are exponentially
decreasing in $|k|$, so using
$\Li_{1-\ell}(z)\approx z$ for small $z$,
and $r_k=e^{-\eps k^2}(1+O(k^3/m^2))$ when $k\leq n^{2/3}$,
we see that doing the substitution for $n^{1/2}\leq k\leq n^{2/3}$
also introduces $O(1)$ error, and that the substitution for $k\geq
n^{2/3}$ gives $o(1)$ error.  Upon multiplying by
$\sqrt{q^3\eps}/(qn)^\ell$, all these errors become $o(1)$.

For the integral approximation \eqref{intapprox}, the integrands in
\eqref{intapprox} are continuous (except at the branch cut) as long as
the curves
$\mt e^{2\pi i j p/q} A e^{2 i \alpha\sqrt{\eps/q}\, k - \eps k^2}$
avoid the point $1$.
Moreover they converge exponentially fast to $0$ when $|k|\to\infty$.
Therefore they are Riemann summable and the error in converting
the sums to integrals tends to zero.
The error introduced by extending the range of integration to the
reals is exponentially small.

What happens when a curve passes through the singularity? 
When $\ell=0$, the integral in expression
\eqref{intapprox} for $\log \Zst$ converges and is an upper bound for
$\log \Zst$: the Riemann sum for $\log \Zst$ 
converges to its integral
on the complement of a small neighborhood of the singularity, and the
Riemann sum near the logarithmic 
singularity has a negligible contribution except
possibly for the point which is closest to the singularity.
If the distance of the closest point to the singularity
is no smaller than $e^{-o(\sqrt{n})}$, then the contribution
of this point is $o(\sqrt{n})$ and so can be ignored.
\end{proof}

\subsection{The distribution of the number of edges of type `$c$'}
\label{sec:distNc}
\enlargethispage*{\baselineskip}

By \lref{ddclz}, 
to first order $\log\Zst$ only depends on the parity of $\tqsp$,
so we define
$$Z_- =\frac12\!\sum_{\substack{\sigma,\tau\\ \text{$\tqsp$ odd}}}\!\!\! \Zst
\ \ \ \ \ \text{and}\ \ \ \ \ 
Z_+ =\frac12\!\sum_{\substack{\sigma,\tau\\ \text{$\tqsp$ even}}}\!\!\! \mst\Zst$$
where $\mst=-1$ if $\sigma=0$ and $\tau=0$, and $\mst=1$ otherwise.
We have $Z=Z_-+Z_+$, and from \pref{2Zsums}, both
$Z_-$ and $Z_+$ have only nonnegative coefficients, so they can
be interpreted as distributions.
\lref{ddclz} shows that typically one of $Z_-$ or $Z_+$ is
exponentially larger than the other one, so the distribution of the
number of $c$ edges is governed by whichever of $Z_-$ or $Z_+$ is dominant.

We use the method of moments to determine the distribution of the
number $N_c$ of type-`$c$' edges.  We saw in
\sref{logZ} how to express the $\ell$th moment of $N_c$
about its mean (call it $C_\ell$)
in terms of $\L_\ell = \ddcl \log Z$, but in \lref{ddclz} we
evaluated $\L_{\ell,\st} = \ddcl \log \Zst$.
Define
$C_{\ell,\st}$ to be the same expression, except with the
$\L_{\ell,\st}$'s replacing the $\L_\ell$'s (see also \eqref{Cstdef} below), 
and similarly define the $\L_{\ell,\pm}$'s and the $C_{\ell,\pm}$'s.

\begin{lemma}
\label{lem:!!}
Under the assumptions of \lref{ddclz}, if the curve
$(-1)^{\tqsp} A^q e^{2 i \alpha x - x^2}$ does not contain the point $1$ and
$\int_{-\infty}^{\infty} \Li_{-1}\big((-1)^{\tqsp} A^q e^{2 i 
\alpha x - x^2}\big) \,dx \neq 0$,
 then $C_{\ell,\st}/C_{2,\st}^{\ell/2}\rightarrow (\ell-1)!!.$
(Here as usual $\ell!!=\ell(\ell-2)\cdots(3)(1)$ when $\ell$ is odd and $\ell!!=0$ when $\ell$ is even.)
\end{lemma}
\begin{proof}
From \lref{ddclz} we have
$\L_{j,\st} = O((qn)^{j}/\sqrt{q^3\eps})$,
and since the above integral is nonzero, 
$C_{2,\st}=\L_{2,\st} = \Theta((qn)^{2}/\sqrt{q^3\eps})$.
Thus each monomial in the polynomial \eqref{Bellpoly} for $C_{\ell,\st}$
has magnitude $O((qn)^\ell/(\sqrt{q^3\eps})^{\text{degree}})$.
Recall that $q^3\eps \ll 1$.  The monomial degree is uniquely maximized by the
$\L_{2,\st}^{\ell/2}$ term ($\ell$ even) or the $\L_{3,\st}
\L_{2,\st}^{(\ell-3)/2}$ term ($\ell$ odd).
Thus when $\ell$ is odd,
$C_{\ell,\st}/C_{2,\st}^{\ell/2} = O((q^3\eps)^{1/4})\rightarrow 0$,
and when $\ell$ is even, $C_{\ell,\st}/C_{2,\st}^{\ell/2}$ tends to
the coefficient of the monomial $\L_{2,\st}^{\ell/2}$ in $C_{\ell,\st}$,
which by \eqref{partcoef} is $(\ell-1)(\ell-3)\cdots (3)(1)$.
\end{proof}
\enlargethispage*{2\baselineskip}

\medskip
\newcommand{\sgn}{\varsigma}
Thus the $C_{\ell,\st}$'s converge to the moments of a Gaussian,
but recall that we cannot view $\Zst$ as a distribution since it may have
some negative coefficients.
Next we show that the moments $C_{\ell,\pm}$ of $Z_\pm$ (which are
genuine distributions) are close to the corresponding $C_{\ell,\st}$'s.
\begin{lemma}
\label{lem:moments}
Suppose $\ell\in\N$, $\sgn=(-1)^{\tqspo}=(-1)^{\tqspt}$,
$(\stoc)\neq(\sttc)$, and either $\sgn=-1$, or else $\sgn=+1$ but $A\geq1$.
Under the assumptions of Lemmas~\ref{lem:ddclz} and \ref{lem:!!},
if $|\L_{1,\sto} - \L_{1,\stt}|\ll n^{5/4}$
then $C_{\ell,\sgn}/C_{2,\sto}^{\ell/2} = (\ell-1)!!+o(1)$.
\end{lemma}
\enlargethispage*{2\baselineskip}
\begin{proof}
Viewing $\Zst$ as a polynomial in $c$,
$\Zst=\sum_{\ii} \coef_{\ii,\st}\, c^\ii$, we have
\begin{equation}\label{Cstdef}C_{\ell,\st}\Zst=
\sum_{\ii} \coef_{\ii,\st}\,(\ii-\L_{1,\st})^\ell c^\ii.
\end{equation}
As $Z_\sgn$ is the average of $\msto Z_{\sto}$ and $\mstt Z_{\stt}$, we have
\begin{align*}
C_{\ell,\sgn} Z_{\sgn} &=\frac12\sum_{\sigma,\tau}\mst \sum_\ii \coef_{\ii,\st}\,c^\ii(\ii-\L_{1,\sgn})^\ell\\
&=\frac12\sum_{\sigma,\tau}\mst\sum_\ii \coef_{\ii,\st}\,c^\ii\sum_{j=0}^\ell
\binom{\ell}{j}(\ii-\L_{1,\st})^j(
\L_{1,\st}-\L_{1,\sgn})^{\ell-j}\\
&=\frac12\sum_{\sigma,\tau}\mst \sum_{j=0}^\ell
\binom{\ell}{j}
C_{j,\st}\Zst(\L_{1,\st}-\L_{1,\sgn})^{\ell-j}.
\end{align*}
Since $\sgn=-1$ or else $\sgn=+1$ but $A\geq1$,
$\msto Z_{\sto}$ and $\mstt Z_{\stt}$ have the same sign,
so $\L_{1,\sgn}$ is a convex combination of $\L_{1,\sto}$ and
$\L_{1,\stt}$, and it too can differ by at most $\ll
n^{5/4} = O(C_{2,\st}^{1/2})$ from them.  Substituting
$C_{j,\st}=((j-1)!!+o(1)) C_{2,\st}^{j/2}$ we get
$$C_{\ell,\sgn} Z_{\sgn} = \frac12\sum_{\sigma,\tau}\mst ((\ell-1)!!+o(1))C_{2,\st}^{\ell/2}\Zst,$$
so
$C_{\ell,\sgn}$ is a convex combination of $((\ell-1)!!+o(1))C_{2,\sto}^{\ell/2}$ and $((\ell-1)!!+o(1))C_{2,\stt}^{\ell/2}$.

Since $C_{2,\sto}=K_{2,\sto}$ and $C_{2,\stt}=K_{2,\stt}$,
\lref{ddclz} and the hypothesis on the $K_2$'s 
implies  $C_{2,\sto}=C_{2,\stt}+O(n^{5/4})$.
Thus $C_{\ell,\sgn}=((\ell-1)!!+o(1))C_{2,\sto}^{\ell/2}.$
\end{proof}

\smallskip
We have not computed $\L_{1,\st}$ to the precision that \lref{moments}
would appear to suggest that we need, but all we really need
is that the two relevant $\L_{1,\st}$'s are quite close.

\begin{lemma}
\label{lem:close}
Suppose $\sgn=(-1)^{\tqspo}=(-1)^{\tqspt}$.
Under the assumptions of \lref{ddclz}, if
$$ -\int_{-\infty}^{\infty} \Li_1(-\sgn A^q e^{2 i \alpha x -x^2}) \,dx
 < -\int_{-\infty}^{\infty} \Li_1(\sgn A^q e^{2 i \alpha x -x^2}) \,dx,
$$
then $|\L_{1,\sto} - \L_{1,\stt}| \leq \exp(-\Theta(\sqrt{n}))$.
\end{lemma}
\begin{proof}
For expository convenience say that the two $\st$'s for which
$(-1)^{\tqsp}=\sgn$ are $01$ and $10$.  From \lref{ddclz}, $Z_{01}$
and $Z_{10}$ dominate $-Z_{00}$ and $Z_{11}$ by a factor of
$\exp(\Theta(\sqrt{n}))$, and then
from \pref{2Zsums}, $|Z_{01}/Z_{10}-1| \leq \exp(-\Theta(\sqrt{n}))$.
Writing $\Zst$ as a polynomial in $c$,
$\Zst=\sum_{\ii} \coef_{\ii,\st}\, c^\ii$, again
from \pref{2Zsums} we have
\begin{align*}
|\coef_{\ii,01}-\coef_{\ii,10}| &\leq -\coef_{\ii,00} + \coef_{\ii,11} \\
|\coef_{\ii,01}\,\ii c^\ii-\coef_{\ii,10}\,\ii c^\ii| &\leq \ii c^\ii (-\coef_{\ii,00} + \coef_{\ii,11}) \leq m n (-\coef_{\ii,00}\,c^\ii  + \coef_{\ii,11}\,c^\ii) \\
\left|\sum_\ii \coef_{\ii,01}\, \ii c^\ii - \sum_\ii \coef_{\ii,10}\, \ii c^\ii \right| &\leq m n \left[-\sum_\ii \coef_{\ii,00}\, c^\ii + \sum_\ii \coef_{\ii,11}\, c^\ii\right] \\
|Z_{01} \L_{1,01} - Z_{10} \L_{1,10}| &\leq m n (-Z_{00}+Z_{11}).
\end{align*}
As $Z_{01}$ and $Z_{10}$ are exponentially close to each other and
exponentially dominate $-Z_{00}$ and $Z_{11}$, it must be that
$\L_{1,01}$ and $\L_{1,10}$ are exponentially close.
\end{proof}

\begin{thm}
\label{thm:gaussian}
When we fix $p$, $q$, $A$, and $\alpha$ while $m\rightarrow\infty$ and
$n\rightarrow\infty$, we have
\begin{equation} \label{p+c}
\log Z = \frac{(m n b c)^{1/4}}{\pi \sqrt{2} (p q)^{3/4}}\left[\max_{\pm} -\int_{-\infty}^{\infty} \Li_{1}\big(\pm A^q e^{2 i \alpha x - x^2}\big) \,dx +o(1)\right] ,
\end{equation}
and with the exceptions noted below,
the number $N_c$ of edges of type `$c$' converges in distribution to a Gaussian,
with mean
\begin{equation} \label{p+c2}
\enc =
 \frac{(m n c)^{3/4}}{\pi \sqrt{2} (p q b)^{1/4}}
\left[
 -\int_{-\infty}^{\infty} \Li_{0}\big(\pm A^q e^{2 i \alpha x - x^2}\big) \,dx  + o(1)\right],
\end{equation}
 and variance
\begin{equation} \label{p+c3}
\vnc =
 \frac{(m n c)^{5/4} (p q)^{1/4}}{\pi \sqrt{2} b^{3/4}}
\left[
 -\int_{-\infty}^{\infty} \Li_{-1}\big(\pm A^q e^{2 i \alpha x - x^2}\big) \,dx + o(1)\right] ,
\end{equation}
where the choice of $\pm$ in \eqref{p+c2} and \eqref{p+c3}
is the value that maximizes \eqref{p+c}.  The exceptions are
\begin{enumerate}
\item When both $+$ and $-$ maximize \eqref{p+c},
the distribution of type-`$c$' edges is a mixture of the two Gaussians
defined above, provided that exceptions 2 and 3 do not also occur.
(In the interest of space we omit the proof about the mixture of Gaussians,
but we can supply it to the interested reader upon request.)
\item If for the dominant choice of $\pm$ in \eqref{p+c}, 
the curve $\pm A^q e^{2 i \alpha x - x^2}$ passes through $1$, we do
not say anything about the distribution of type-`$c$' edges, but the
formula for $\log Z$ remains valid.  (We have reason to believe, but
have not proved, that this scenario never occurs.)
\item In the event that the integral in \eqref{p+c3} evaluates to $0$,
the formulas are still valid, but we no longer claim that the distribution
is a Gaussian.  (We believe that this scenario never occurs, \lref{var>}
in \sref{>'s} rules it out when $A\leq1$.)
\end{enumerate}
\end{thm}
(See \aref{notation} for a review and discussion of the parameters $A$ and $\alpha$.)

\begin{proof}
Immediate from \lref{ddclz}, \lref{!!},
\lref{moments} (with the fact that \lref{-} 
shows the dominant choice of $\pm$ to be $-$ when $A\leq 1$), \lref{close},
the fact that the moments are those of a
Gaussian random variable, the method of moments,
approximation~\eqref{q3eps} for $1/\sqrt{q^3\eps}$,
and $qn\approx \sqrt{p q m n c/b}$.

In the event that one of the curves $\pm A^q e^{2 i \alpha x - x^2}$
passes through the point $1$, we still obtain \eqref{p+c} from \lref{ddclz} because each
such curve has two corresponding $\Zst$'s, at most one of which can have a multiplicand closer than $e^{-o(\sqrt{n})}$ to $0$.
\end{proof}

\tref{intro} follows by plugging into \tref{gaussian}
$A=1$ and $\alpha=0$, using the explicit evaluation of the integrals in
\sref{Zxint}, and using $-\Li_{3/2}(-1)>-\Li_{3/2}(1)$,
$-\Li_{-1/2}(-1)>0$, and $\Li_\nu(-1)=(2^{1-\nu}-1)\zeta(\nu)$.

\begin{figure}[htbp]
\centerline{\epsfig{figure=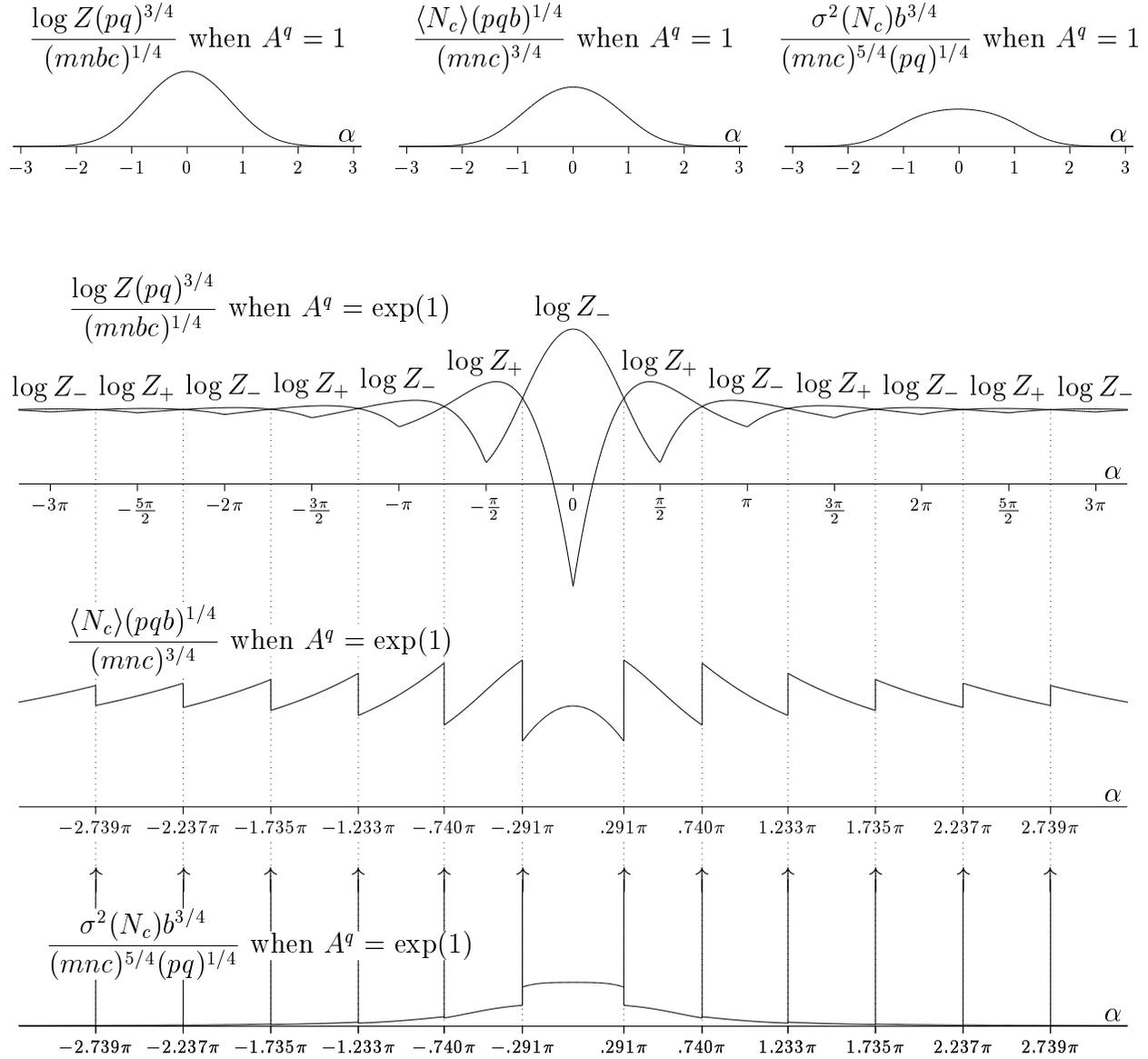}}
\caption{
Anatomy of the resonant spikes.  In the upper panel we show the curves
for $\log Z$, $\enc$, and $\vnc$ as a function of $\alpha$ when $A=1$.
(The parameters $A$ and $\alpha$ are reviewed in \aref{notation}.)
When $A\leq1$, $Z_-$ is dominant,
and all three curves are analytic and unimodal.
The situation is quite different when $A>1$.
In the next panel we show the curves for $\log Z_+$ and $\log Z_-$
as a function of $\alpha$ when $A^q=\exp(1)$ ($\log Z$ is the max of these
two curves).  When $A>1$ there are singularities
in $\log Z_+$ and $\log Z_-$ that occur when their corresponding spirals
cross the singularity, that is, when $\alpha=(\pi/\sqrt{\log A^q}) \Z$ for
$\log Z_+$ and when $\alpha=(\pi/\sqrt{\log A^q}) (\Z+1/2)$ for $\log Z_-$.
In the lower two panels we show the curves for $\enc$ and $\vnc$.
These curves were computed using \eqref{p+c2} and \eqref{p+c3} as explicitly
evaluated in \eqref{IntLi}
with whichever of $Z_+$ or $Z_-$ is significant.
The ``crossover points'',
where $\log Z_-$ and $\log Z_+$ alternate in significance,
are in a sense a phase transition within a phase transition.
At each crossover point, the curve for $\log Z$ is continuous but nonanalytic,
the curve for $\enc$ is discontinuous,
and the curve for $\vnc$ has a delta function.
When a spiral hits the singularity, the corresponding $Z_\pm$ appears to be
the insignificant one. 
The large-$\alpha$ asymptotics of the curves for $\log Z$, $\enc$, and $\vnc$
are given in \tref{Z,alpha=infinity}.}
\label{fig:Zsig}
\end{figure}

\section{Understanding the resonant spikes}
\label{sec:spike}

The resonant spikes, an example of which is shown in \fref{Zsig},
exhibit nontrivial behavior.  We saw already that this behavior is
determined by the integrals $$\int_{-\infty}^{\infty}
\Li_\nu\big(\beta e^{2 i\alpha x -x^2}\big) \,dx,$$ with $\nu=1$
governing $\log Z$, $\nu=0$ governing $\enc$, and
$\nu=-1$ governing $\vnc$.  We start by explicitly evaluating
these integrals in \sref{Zxint}, and then investigate some of their
properties in \sref{>'s}, \sref{bigalpha}, and \sref{crossover}.
These subsections may be read in any order.

\subsection{The integral $\displaystyle\int_{-\infty}^{\infty} \Li_\nu\big(\beta e^{2 i\alpha x -x^2}\big) \,dx$}\label{sec:Zxint}

We assume that $\beta\in\C$, $\alpha\in\R$, and $\nu\in\C$, though later we restrict $\nu$ to integers $\leq 1$.
For convenience we assume $\alpha\geq0$ since the integral is an even function of $\alpha$.
When $|\beta|>1$ there is a branch cut ($\Li_\nu(z)$ has a branch
cut $\{z\in\R: z\geq 1\}$) that may be encountered when we
vary $x$, so we need to specify which branch of the polylogarithm we
are integrating over.  Since the principal branch is the one for which
$\Li_\nu\big(\beta e^{2i\alpha x - x^2}\big)\rightarrow 0$ as
$x\rightarrow\pm\infty$, we specify that the integrand is the
principal branch of $\Li_\nu$, even though this may make the
integrand only piecewise analytic as a function of $x$ as it ranges
from $-\infty$ to $\infty$.  In the event that $\alpha=0$ and $\beta$ is
real and $\geq 1$, so as to ensure continuity in $\alpha$, we specify that
$\beta e^{2i\alpha x - x^2}$ lies above the branch cut for positive $x$
and below the cut for negative $x$.

We set $z=x+iy$ ($x,y\in\R$) and integrate instead within the complex plane.
Rather than integrate $\Li_\nu\big(\beta e^{2 i\alpha z -z^2}\big)$ along the
real axis $\Im z=0$, it is more convenient to the integrate along the
line $\Im z = \alpha$.  When we deform the contour of integration and
set $z=x+i\alpha$ we get
$$ \int_{-\infty}^\infty \Li_\nu\big(\beta e^{2i\alpha x-x^2}\big) \,dx =
   \int_{-\infty}^\infty \Li_\nu\big(\beta e^{-\alpha^2} e^{-x^2}\big) \,dx +
   \parbox{2.25in}{terms from the singularities and branch cuts of $\Li_\nu\big(\beta e^{2i\alpha z-z^2}\big) $ in the complex $z$-plane.}$$

In particular if $|\beta|\leq 1$ there are no singularities or branch
cuts encountered when deforming the contour of integration, so there
are no additional terms.  We shall evaluate the main term in first.
For now assume that $|\beta| e^{-\alpha^2}<1$ so that the series
expansion of the polylogarithm is absolutely convergent.  This enables
us to write
\begin{align*}
\int_{-\infty}^{\infty} \Li_\nu\big(\beta e^{-\alpha^2} e^{-x^2}\big) \,dx
 &= \int_{-\infty}^{\infty}
     \sum_{n=1}^{\infty}
     \frac{(\beta e^{-\alpha^2})^n}{n^\nu}
     e^{-n x^2} \,dx \\
 &= \sum_{n=1}^{\infty}
     \frac{(\beta e^{-\alpha^2})^n}{n^\nu}
     \int_{-\infty}^{\infty}
     e^{-n x^2} \,dx \\
 &= \sum_{n=1}^{\infty}
     \frac{(\beta e^{-\alpha^2})^n}{n^\nu}
     \frac{\sqrt{\pi}}{\sqrt{n}} \\
 &= \sqrt{\pi} \Li_{\nu+1/2}\big(\beta e^{-\alpha^2}\big).
\end{align*}

The singularities of $\Li_\nu\big(\beta e^{2i\alpha z-z^2}\big)$ in the complex $z$-plane occur when
\begin{align*}
\beta e^{2i\alpha z-z^2} &= 1\\
2i\alpha z-z^2 &= -\log\beta-2\pi i k\\
z &= i\alpha\pm i\sqrt{\alpha^2-\log\beta-2\pi i k}
\end{align*}
where $k$ is an integer.
Since we are moving the contour between $\Im z =0$ and $\Im z =
\alpha$ (see \fref{branch}),
\begin{figure}[htbp]
\centerline{\epsfig{figure=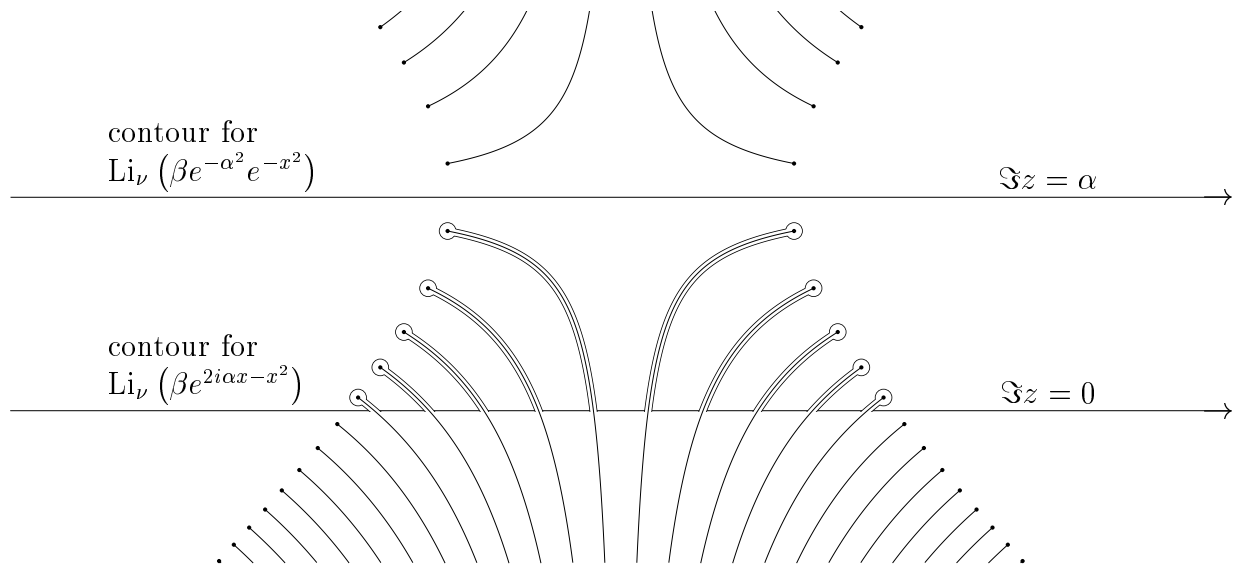}}
\caption{Singularities, branch cuts, and contours of integration for $\Li_\nu\big(\beta e^{2 i\alpha z -z^2}\big)$.  When pushing the contour from $\Im z=\alpha$ to $\Im z=0$, the contour must deform to travel around the singularities and branch cuts, so these contribute to the integral.}
\label{fig:branch}
\end{figure}
the relevant singularities are of the form $i\alpha -
i\sqrt{\alpha^2-\log\beta-2\pi i k}$ where the principal square root
is taken.  (If $\arg\beta=0$, $\alpha^2<\log|\beta|$, and $k=0$ then both roots are relevant.)
The branch cuts of the integrand occur where $\beta
e^{2i\alpha z-z^2}$ is real and $\geq 1$.
With $z=x+i y$ and $x,y\in\R$ these are
\begin{align*}
2\alpha x-2x y &= -\arg\beta-2\pi k & -2\alpha y -x^2+y^2 &\geq-\log|\beta| \\
x &= \frac{1}{2}\frac{2\pi k+\arg\beta}{y-\alpha} & (y-\alpha)^2-x^2&\geq\alpha^2-\log|\beta|
\end{align*}

To identify if and where the $k$th cut intersects the line $\Im z=0$
we set $y=0$ and find that there is an intersection at
$$-\frac{2\pi k+\arg\beta}{2\alpha}$$ provided
$$ (2\pi k+\arg\beta)^2 \leq 4\alpha^2\log|\beta|.$$

For the moment suppose that $|\beta| e^{-\alpha^2}<1$ so that there are
no singularities or branch cuts of the integrand on the line $\Im
z=\alpha$.  Push this contour downwards toward the line $\Im z =0$,
except that when a branch cut is encountered, the contour must go
around the branch cut, and travels along the cut upwards on its left
side in downwards on its right side (\fref{branch}).
The integrand increases by
$(2\pi i/\Gamma(\nu)) (\log(\beta e^{2i\alpha z-z^2}))^{\nu-1}$ when
going from the left side of the cut to the right side, so the net
contribution of this cut to the integral is $$\frac{2\pi i}{\Gamma(\nu)}
\int_C (\log(\beta e^{2i\alpha z-z^2}))^{\nu-1} \,dz $$ (plus another
term from the singularity) where the contour $C$ is a branch cut
travelling from the branch point to the point where the branch cut
intersects the line $\Im z = 0$.  For the $k$th branch cut the
endpoints of integration are $i\alpha-i\sqrt{\alpha^2-\log\beta-2\pi i
k}$ and $-(2\pi k +\arg\beta)/(2\alpha)$.  When $\nu=1$
the singularity does not contribute and this integral
is easy to evaluate and we get $$2\pi i\left[
-i\alpha+i\sqrt{\alpha^2-\log\beta-2\pi i k} -\frac{2\pi k+\arg\beta}{2\alpha}
\right].$$
Since this term gets added to the integral over $\Im z=0$, when we evaluate the integral over $\Im z=0$ we subtract this term from the integral over $\Im z=\alpha$.
Thus
\begin{equation}\label{5}
\begin{split}
\int_{-\infty}^\infty \Li_1\big(\beta e^{2i\alpha x-x^2}\big) \,dx =
 \sqrt\pi &\Li_{3/2}\big(\beta e^{-\alpha^2}\big) - \\ &2\pi i
\!\!\!\!\!\!\!\!\!\!\!\!\!\!\!\!
\sum_{\substack{k\in\Z\\(2\pi k+\arg\beta)^2 \leq 4\alpha^2\log|\beta|}}
\!\!\!\!\!\!\!\!\!\!\!\!\!\!\!\!
\left[
-i\alpha+i\sqrt{\alpha^2-\log\beta-2\pi i k} -\frac{2\pi
k+\arg\beta}{2\alpha}\right].
\end{split}
\end{equation}
There is a term in the summation for each time the spiral $\beta e^{2i\alpha
x-x^2}$ encloses the singularity at $1$.
This formula is valid when $|\beta| e^{-\alpha^2}<1$,
and we would like to show that it is valid (for real $\beta$ and $\alpha$)
without this restriction.
We can do this via analytic continuation in (complex) $\beta$ and $\alpha$.
While both sides of equation \eqref{5} are analytic in $\alpha$
(except where the spiral
$\beta e^{2i\alpha x-x^2}$ crosses the singularity),
neither side of \eqref{5} is analytic in $\beta$.
The $\arg\beta$ terms on the right-hand side are nonanalytic, and
the left-hand side is nonanalytic due to the branch cut in $\Li_1(z)$ ---
as $\arg\beta$ is varied, portions of the integrand cross the branch cut.
To remedy this problem, we deform the branch cut as $\arg\beta$
changes, so that the cut consists of the arc $\{z\in\C:\text{$|z|=1$
and $0\leq\arg z\leq\arg\beta$ or $0\geq\arg z\geq\arg\beta$}\}$
together with the ray $\{z\in\C:\text{$|z|\geq 1$ and $\arg
z=\arg\beta$}\}$.
So long as the spiral avoids the singularity, for any value of $x$,
the integrand does not cross the moving branch as $\beta$ and $\alpha$
are changed.
This re-interpreted integral (with the moving branch cut) is then
analytic except where
the spiral $\beta e^{2i\alpha x-x^2}$ crosses the singularity,
and the $\arg\beta$ terms on the right-hand side become constant.
The $\sqrt\pi \Li_{3/2}\!\big(\beta e^{-\alpha^2}\big)\!$ term on
the right-hand-side of \eqref{5} is nonanalytic when
$\beta e^{-\alpha^2}=1$, but this nonanalyticity is cancelled by
the nonanalyticity in the $k=0$ term in the summation
(by \eqref{Li-pole} in \aref{polylogs}, which reviews polylogarithms).
Thus the right-hand-side of the re-interpreted \eqref{5} is also analytic
whenever the spiral $\beta e^{2i\alpha x-x^2}$ avoids the singularity.
For a given $\beta$ and $\alpha$, we may continuously decrease
$|\beta|$ and increase $\alpha$
until $|\beta|e^{-\alpha^2}<1$ (where we already know that \eqref{5} is
valid), while keeping the spiral from crossing the singularity.
Thus \eqref{5} is valid for all complex $\beta$ and $\alpha$ for which
the spiral avoids the singularity.  As both sides of \eqref{5} are
continuous, \eqref{5} is also valid when the spiral contains the singularity.
We revert to the principal branch cut interpretation of the integral,
restoring the $\arg\beta$ terms on the right-hand side of \eqref{5}.

When $\beta$ is real, upon summing over $k$, the terms
$(2\pi k+\arg\beta)/(2\alpha)$ cancel: when $\beta<0$ the $k$th term
cancels the $(-1-k)$th term, and when $\beta>0$ the $k$th term cancels
the $-k$th term, while the $0$th term does not contribute.
For real $\beta$ we get
\begin{equation}
\label{IntLi1}
\int_{-\infty}^{\infty} \Li_1\big(\beta e^{2 i\alpha x -x^2}\big) \,dx
 = \sqrt{\pi} \Li_{3/2}\big(\beta e^{-\alpha^2}\big) -
2\pi
\!\!\!\!\!\!\!\!\!\!\!\!\!\!\!\!
\sum_{\substack{k\in\Z\\(2\pi k+\arg\beta)^2 \leq 4\alpha^2\log|\beta|}}
\!\!\!\!\!\!\!\!\!\!\!\!\!\!\!\!
     \left(\alpha-\sqrt{\alpha^2-\log\beta-2\pi i k}\right).
\end{equation}
When $\beta e^{-\alpha^2}>1$ the $k=0$ ``branch-cut term'' exactly
cancels the imaginary part of the $\Li_{3/2}$ term.
Applying $\beta \frac{\partial}{\partial\beta}$ multiple times to
\eqref{IntLi1} we find for real $\beta$ and nonpositive integer $\nu$
\begin{equation}
\label{IntLi}
\int_{-\infty}^{\infty} \Li_\nu\!\big(\beta e^{2 i\alpha x -x^2}\big)\! \,dx
 = \sqrt{\pi} \Li_{\nu+1/2}\!\big(\beta e^{-\alpha^2}\big)\!
     -\sqrt{\pi} \Gamma({\textstyle \frac12\!-\!\nu}) 
 \sum_{\substack{k\in\Z\\ \makebox[0pt]{$\scriptstyle(2\pi k+\arg\beta)^2 \leq 4\alpha^2\log|\beta|$}}}
[\alpha^2-\log\beta+2\pi i k]^{-1/2+\nu}.
\end{equation}

\subsection{Inequalities}
\label{sec:>'s}

To use \tref{gaussian} we need to know which of
$Z_{\pm}$ is significant.
We already saw that they can alternate in significance (as $\alpha$ is
varied) when $A>1$.  We claim that when $A\leq 1$, or else $A>1$ but
$\alpha=0$, that $Z_-$ is the significant one.
\begin{lemma}
\label{lem:-}
If $0<\beta\leq 1$, or else $\beta>1$ but $\alpha=0$, then
\begin{equation}
\label{z-}
   -\int_{-\infty}^{\infty} \Li_1\big(\beta e^{2 i \alpha x -x^2}\big) \,dx
< -\int_{-\infty}^{\infty} \Li_1\big(-\beta e^{2 i \alpha x -x^2}\big) \,dx.
\end{equation}
\end{lemma}
\begin{proofm}
Under either of the hypotheses, equation~\eqref{IntLi1} simplifies to
$$
\int_{-\infty}^{\infty} \Li_1\big(\pm\beta e^{2 i\alpha x -x^2}\big) \,dx
 = \sqrt{\pi} \Li_{3/2}\big(\pm\beta e^{-\alpha^2}\big)
 = \int_{-\infty}^{\infty} \Li_1\big(\pm\beta e^{-\alpha^2-x^2}\big) \,dx,
$$
but
\begin{equation}
 -\Re\Li_1\!\big(\beta e^{-\alpha^2-x^2}\big)\!
 = \Re\log\!\big(1-\beta e^{-\alpha^2-x^2}\big)\!
 <
 \Re\log\!\big(1+\beta e^{-\alpha^2-x^2}\big)\!
 = -\Re\Li_1\!\big(\!-\beta e^{-\alpha^2-x^2}\big). \tag*{\qed}
\end{equation}
\end{proofm}

To show that the distribution of $N_c$ is a Gaussian, we need to know
that the integral expression in \eqref{p+c3} (for the dominant choice
of $\pm$) is nonzero.  The integral clearly cannot be negative,
because it has an interpretation in terms of variance, but \textit{a
priori\/} it could be zero, in which case the $o(1)$ error term would
control the variance, and we would be unable to characterize the
distribution of $N_c$.  Ideally we would like to show that it always
positive, but we only do this for $A\leq 1$, or else $A>1$ but $\alpha=0$.
Recall that under these conditions $Z_-$ is dominant.
\begin{lemma}
\label{lem:var>}
If $0<\beta\leq 1$, or else $\beta>1$ but $\alpha=0$, then
$$ -\int_{-\infty}^{\infty} \Li_{-1}\big(-\beta e^{2 i \alpha x - x^2}\big) \,dx > 0.$$
\end{lemma}
\begin{proof}
As above,
under either of the hypotheses, equation~\eqref{IntLi} simplifies to
$$
\int_{-\infty}^{\infty} \Li_{-1}\big(-\beta e^{2 i\alpha x -x^2}\big) \,dx
 = \sqrt{\pi} \Li_{-1/2}\big(-\beta e^{-\alpha^2}\big)
 = \int_{-\infty}^{\infty} \Li_{-1}\big(-\beta e^{-\alpha^2-x^2}\big) \,dx,
$$
but
$
 -\Li_{-1}\big(-\beta e^{-\alpha^2-x^2}\big)
 = {\beta e^{-\alpha^2-x^2}}/{\big(1+\beta e^{-\alpha^2-x^2}\big)^2}
 >0.
$
\end{proof}

\begin{lemma}
The curve for $\log Z$ attains its (unique) maximum value when $\alpha=0$.
\end{lemma}
\begin{proof}
When $0<\beta$ we have $\log(1+\beta e^{-x^2}\big) \geq \Re \log(1\pm \beta e^{2 i \alpha x-x^2}\big)$, with strict inequality when $\pm e^{2 i \alpha x}\neq 1$, so
$
 -\int_{-\infty}^{\infty} \Li_1\big(\mp\beta e^{2 i \alpha x -x^2}\big) \,dx
\leq -\int_{-\infty}^{\infty} \Li_1\big(-\beta e^{-x^2}\big) \,dx
$,
with strict inequality unless both $\alpha=0$ and $\mp=-$.
\end{proof}
\medskip

\noindent
In contrast, the curve for $\enc$ does not always attain its
maximum value when $\alpha=0$.

\subsection{Asymptotics for large $\alpha$}
\label{sec:bigalpha}

In \fref{Zsig} the curves for $\log Z$ and $\enc$ appear
to asymptote out to a positive constant, while the curve for $\vnc$
decays to $0$.  The following theorem gives the large-$\alpha$
asymptotics of these curves which are given by the integrals in
\eqref{p+c}, \eqref{p+c2}, and \eqref{p+c3} in \tref{gaussian}.

\newcommand{\ab}{\beta}
\begin{thm}
\label{thm:Z,alpha=infinity}
For positive $\beta$ and real $\alpha$,
\begin{align*}
\lim_{|\alpha|\rightarrow\infty}
-\int_{-\infty}^{\infty} \Li_1\big(\pm\beta e^{2 i \alpha x - x^2}\big) \,dx
&= \begin{cases} \frac{4}{3} (\log\ab)^{3/2} & \ab\geq 1 \\
  0 & \ab\leq 1\end{cases} \\
\lim_{|\alpha|\rightarrow\infty}
-\int_{-\infty}^{\infty} \Li_0\big(\pm\beta e^{2 i \alpha x - x^2}\big) \,dx
&= \begin{cases} 2 (\log\ab)^{1/2} & \ab\geq 1 \\
  0 & \ab\leq 1\end{cases}\\
\lim_{|\alpha|\rightarrow\infty}
-\alpha^2 \int_{-\infty}^{\infty} \Li_{-1}\big(\pm\beta e^{2 i \alpha x - x^2}\big) \,dx
&= \begin{cases} (\log\ab)^{1/2} & \ab\geq 1 \\
  0 & \ab\leq 1\end{cases}
\end{align*}
\end{thm}
\begin{proof}
These integrals are evaluated in \eqref{IntLi1} and \eqref{IntLi}.
When $|\alpha| \rightarrow\infty$, the main terms in \eqref{IntLi1}
and \eqref{IntLi}, $-\sqrt{\pi} \Li_{\nu+1/2}(\pm \ab e^{-\alpha^2})$
for $\nu=1,0,-1$, become negligible.  Thus we only have to compute the
``branch-cut terms'' when $\beta>1$.  For \eqref{IntLi1} these terms
can be approximated for large $\alpha$ as follows.
Let $\gamma=\log \ab - \alpha^2$.  The identity
$$\sqrt{x+iy}+\sqrt{x-iy}=\sqrt{2\big(x+\sqrt{x^2+y^2}\big)}$$
is easily verified by squaring both sides, and allows us to rewrite
the sum in \eqref{IntLi1} as
$$\sum_{|k|<|\alpha|\sqrt{\log \ab}/\pi}2\pi\left(|\alpha|-\sqrt{
\big(-\gamma+\sqrt{\gamma^2+4\pi^2k^2}\big)/2}\right)$$
where $k$ is integer for positive $\beta$ (for $Z_+$) or
half-integer for negative $\beta$ (for $Z_-$).

Suppose $|\alpha| \gg \sqrt{\log \ab}$.  Then $k$, which ranges up to
$(|\alpha|/\pi) \sqrt{\log \ab}$, is much smaller than $|\gamma| =
(1+o(1))\alpha^2$, and we have
\begin{align*}
  \sqrt{\big(-\gamma+\sqrt{\gamma^2+4 \pi^2 k^2}\big)/2} &=
  \sqrt{-\gamma} \sqrt{\big(1+\sqrt{1+4 \pi^2 k^2/\gamma^2}\big)/2} \\
 &=
  \sqrt{-\gamma} \sqrt{1 + (\pi^2+o(1)) k^2/\gamma^2} \\
 &=
  \sqrt{\alpha^2-\log \ab} (1 + (\pi^2/2+o(1)) k^2/\gamma^2) \\
 &=
  |\alpha| (1-(1/2+o(1))\log \ab/\alpha^2) (1 + (\pi^2/2+o(1)) k^2/\alpha^4) \\
 &=
  |\alpha| - \frac{1+o(1)}{2 |\alpha|} (\log \ab - \pi^2 k^2/\alpha^2).
\end{align*}
In particular
\begin{align*}
   2 \pi \left(|\alpha|-\sqrt{\big(-\gamma+\sqrt{\gamma^2+4 \pi^2 k^2}\big)/2}\right) &=
   \frac{\pi+o(1)}{|\alpha|} (\log \ab - \pi^2 k^2/\alpha^2)
\end{align*}
and we sum this over $k$'s for which the quantity is positive.
If $|\alpha|\gg 1/\sqrt{\log\ab}$ then we can approximate this
sum with an integral, and the integral is
$$\frac{\pi+o(1)}{|\alpha|} \frac{2}{3} \log \ab
\times \frac{2|\alpha|}{\pi}\sqrt{\log \ab} = (1+o(1))\frac{4}{3} (\log \ab)^{3/2}.$$

These are the asymptotics for $\nu=1$.  For integer $\nu\leq 0$ the
range of $k$ is the same as it is for $\nu=1$, i.e.\ up to about
$|\alpha| \sqrt{\log \ab}/\pi$.  Thus we need to evaluate
$$
\sqrt{\pi} \Gamma(1/2-\nu) \sum_{|k|<|\alpha| \sqrt{\log \ab}/\pi} [\alpha^2-\log \ab+2\pi i k]^{-1/2+\nu}.
$$
For $\alpha\gg\sqrt{\log\ab}$, $k\ll\alpha^2$ so that $$
[\alpha^2-\log \ab+2\pi i k]^{-1/2+\nu} + [\alpha^2-\log \ab-2\pi i k]^{-1/2+\nu} \approx 2 |\alpha|^{2\nu-1}$$
and if $\alpha\gg 1/\sqrt{\log\ab}$ the summation is approximately an integral which is asymptotically
$$ \frac{2}{\sqrt{\pi}} \Gamma(1/2-\nu) \alpha^{2\nu} \sqrt{\log \ab}.$$
Taking $\nu=0$ and $\nu=-1$ give the desired results for the $\enc$ and $\vnc$ integrals.
\end{proof}
\medskip

\subsection{Crossover locations}
\label{sec:crossover}

If $A\leq 1$ then $Z_-$ always exceeds $Z_+$.  But for
larger $A$ there are crossover values for $\alpha$ at which
$Z_-$ and $Z_+$ alternate in significance.
Since $Z_-$ and $Z_+$ count configurations in different $\Z_2$-homology classes
(see \eqref{Ztable}), we conclude that
the crossover values are the places where the typical homology class
of a lattice path changes.  So there is a phase transition
at these points: the topology of a typical configuration changes.

The crossover values $\alpha$ for a fixed $\beta=A^q$
satisfy an implicit (and transcendental) equation.
Instead of solving these equations directly
we can derive an analytic expression for
the crossover $\alpha$'s as a function of $\gamma =
\log \beta-\alpha^2$, and then given $\gamma$ and the crossover $\alpha$
we can calculate the corresponding $\beta$.  In this way we can
parametrically plot these critical pairs $(\beta,\alpha)$ as a
function of $\gamma$.
For example, for the $0$th crossover we have
\begin{align*}
  \eps^{1/2}q^{3/2}Z_- +o(1)&= -\sqrt{\pi}\Li_{3/2}(-\beta e^{-\alpha^2})\\
 \eps^{1/2}q^{3/2}Z_+ +o(1)&= -\sqrt{\pi}\Li_{3/2}( \beta e^{-\alpha^2})
 +2\pi(\alpha-\sqrt{\alpha^2-\log \beta})\\
\intertext{from which we can solve}
\alpha &= \frac{\Li_{3/2}(e^\gamma) - \Li_{3/2}(-e^\gamma)}{2\sqrt{\pi}} + \sqrt{-\gamma}.\\
\intertext{Similarly, for the $r$th crossover we have}
(-1)^r \alpha &= \frac{\Li_{3/2}(e^\gamma) - \Li_{3/2}(-e^\gamma)}{2\sqrt{\pi}}
      +\sum_{k=-r}^r (-1)^k \sqrt{-\gamma+k\pi i}.
\end{align*}

Next let us approximate these crossovers for $r$ fixed and $\gamma$
large.  For large $\gamma$ we can substitute $\nu=3/2$ into the asymptotic
series expansions \eqref{Li(-x)} and \eqref{Li(x)} for $\Li_\nu(\pm
e^\gamma)$ to write
\begin{equation*}
       -\Li_{3/2}(-e^{\gamma}) \approx \frac{4}{3}  /\pi^{1/2} \gamma^{3/2}
                     +\frac{1}{6}   \pi^{3/2}/\gamma^{1/2}
                     +\frac{7}{480} \pi^{7/2}/\gamma^{5/2}
                     +\cdots
\end{equation*}
\begin{equation*}
        \Li_{3/2}(e^{\gamma}) \approx -\frac{4}{3}  /\pi^{1/2} \gamma^{3/2}
                     -2 \sqrt{\pi} i \sqrt{\gamma}
                     +\frac{1}{3}   \pi^{3/2}/\gamma^{1/2}
                     +\frac{1}{60}  \pi^{7/2}/\gamma^{5/2}
                     +\cdots
\end{equation*}
so that
\begin{align*}
(-1)^r \alpha &\approx \frac{\pi}{4}\gamma^{-1/2} +
      \sum_{k=1}^r (-1)^k \left[\sqrt{-\gamma+k\pi i}+\sqrt{-\gamma-k\pi i}\right] \\
 &= \frac{\pi}{4} \gamma^{-1/2} +
      \sum_{k=1}^r (-1)^k \sqrt{2\big(-\gamma+\sqrt{\gamma^2+k^2\pi^2}\big)} \\
 &\approx \frac{\pi}{4} \gamma^{-1/2} +
      \sum_{k=1}^r (-1)^k \sqrt{2(-\gamma+\gamma(1+k^2\pi^2/(2\gamma^2)))} \\
 &= \frac{\pi}{4} \gamma^{-1/2} +
      \sum_{k=1}^r (-1)^k \sqrt{k^2\pi^2/\gamma} \\
 &= \frac{\pi}{4} \gamma^{-1/2} +
      \sum_{k=1}^r (-1)^k k\pi/\gamma^{1/2} \\
 &= (-1)^r (r/2+1/4) \pi \gamma^{-1/2}.
\end{align*}
Since $\alpha$ is small, $\gamma\approx \log \beta$, and so the $r$th crossover
occurs at $\alpha \approx (r/2+1/4) \pi / \sqrt{\log \beta}$.

By comparison the nonanalyticities in the curve for $Z_-$ occur exactly at
$\alpha=\pi/\sqrt{\log \beta} (\Z+1/2)$, and the nonanalyticities in
the curve for $Z_+$ occur exactly at $\alpha=\pi/\sqrt{\log \beta} \Z$.

\section{Open problems}\label{sec:open}

Our analysis for of the lattice paths at the critical point is geared
to regions whose aspect ratio is close to a simple rational number,
but we do not know how to treat ``irrational domains'' which do not have
a simple rational approximation.  Consider for
instance the case $a=1$, $b=c=1/2$ when the side lengths $m$ and $n$
are successive Fibonacci numbers.  Then the aspect ratio of the region
is very close to the Golden ratio, which is not well approximated by
simple rationals.  In this case we believe that the partition function
$Z$ is $\Theta(1)$ (about $2.1$), so that with $\Theta(1)$ probability
there are no lattice paths.  In the event that there are lattice
paths, we believe that they connect up into $\Theta(1)$ loops each of
length $\Theta(n^{4/3})$.  We have also found empirically that as one
varies the aspect ratio of large regions, the smallest value that the
partition function $Z$ takes on is very close to $2$.  We do not
know how to prove any of these conjectures.

In \fref{spike} we plotted $\log Z$ and $\enc$ versus $n/m$ when
$a=1$, $b=c=1/2$.  As predicted there are spikes in $\log Z$ and
$\enc$ when the aspect ratio $n/m$ is near a simple rational $p/q$.
But there also appear to be flanking secondary spikes in $\enc$ when
$n/m \approx p/q$ but $|n/m-p/q|\gg 1/\sqrt{n}$.  We do not
understand this phenomenon.

When $a=1$,
$b=c=1/2+\Theta(1/n)$ and the aspect ratio is nearly a simple rational
$p/q$, the partition function and edge density as a function of $\alpha$ are
nonanalytic at certain points.  We know that the homology type of the
strands changes at these nonanalyticities.  We conjecture that when
$\alpha=0$ each loop winds around exactly $p$ times horizontally and $q$
times vertically.  For large enough values of $\alpha$, when one follows a
loop for $m p + n q$ steps, one does not return to the starting point.
We conjecture that the number of nonanalyticities between $\alpha$ and $0$
determines how many strands away one ends up after following a loop
for $m p + n q$ steps.

Despite our explicit formulas for $\log Z_{\pm}$ and for the expected
value and variance of the number of edges, there are some basic
properties about these functions that we have not been able to derive.
For example, we conjecture that when our formula for $\log Z_-$ is
nonanalytic, our formula for $\log Z_+$ gives a \textit{strictly\/}
larger value, and vice
versa, and that the formula for $\vnc$ is strictly positive.
We also conjecture that for $A>1$ our formula for $\log Z$ as
a function of $\alpha$ always exceeds its limiting value as
$\alpha\rightarrow\infty$.  These conjectures seem fairly evident
from the graphs in \fref{Zsig}, but it is not obvious how to
prove them.  One would also like to determine for fixed $A>1$ the
maximum and minimum values given by our formula for $\enc$,
since both of these can be different than the
limiting $\alpha\rightarrow\infty$ value.

\bibliography{kw}

\begin{thebibliography}{10}

\bibitem{bateman-erdelyi:vol1}
Harry Bateman and A.~Erd\'elyi \textit{et al.}
\newblock {\em Higher Transcendental Functions}, volume~1.
\newblock McGraw-Hill, 1953.
\newblock Based in part on notes left by Bateman, edited by Erd\'elyi
  \textit{et al}.

\bibitem{blote-hilhorst:roughening}
H.~W.~J. Bl{\"o}te and H.~J. Hilhorst.
\newblock Roughening transitions and the zero-temperature triangular {I}sing
  antiferromagnet.
\newblock {\em Journal of Physics \textrm{A}}, 15(11):L631--L637, 1982.

\bibitem{BCK}
C.~Borgs, J.~T. Chayes, and C.~King.
\newblock Meissner phase for a model of oriented flux lines.
\newblock {\em Journal of Physics \textrm{A}}, 28(23):6483--6499, 1995.

\bibitem{cerf-kenyon:wulff}
Rapha{\"e}l Cerf and Richard Kenyon.
\newblock The low-temperature expansion of the {Wulff} crystal in the {3D}
  {Ising} model.
\newblock {\em Communications in Mathematical Physics}, 222:147--179, 2001.

\bibitem{chang-peres:zeta}
Joseph~T. Chang and Yuval Peres.
\newblock Ladder heights, {Gaussian} random walks and the {Riemann} zeta
  function.
\newblock {\em Annals of Probability}, 25(2):787--802, 1997.

\bibitem{cohn-kenyon-propp:variational}
Henry Cohn, Richard Kenyon, and James Propp.
\newblock A variational principle for domino tilings.
\newblock {\em Journal of the American Mathematical Society}, 14(2):297--346,
  2001.
\newblock arXiv:math.CO/0008220.

\bibitem{Comtet}
Louis Comtet.
\newblock {\em Advanced Combinatorics}.
\newblock D. Reidel Publishing Co., 1974.

\bibitem{dN}
Marcel den Nijs.
\newblock The domain wall theory of two-dimensional commensurate-incommensurate
  phase transitions.
\newblock In C.~Domb and J.~L. Lebowitz, editors, {\em Phase Transitions and
  Critical Phenomena}, volume~12, pages 219--333. Academic Press, 1988.

\bibitem{dingle:bose-einstein}
R.~B. Dingle.
\newblock On the {Bose-Einstein} integrals {${\cal B}_p(\eta) = (p!)^{-1}
  \int_0^\infty \varepsilon^p (e^{\varepsilon-\eta}-1)^{-1} d\varepsilon$}.
\newblock {\em Applied Scientific Research}, 6(4):240--244, 1957.

\bibitem{dingle:fermi-dirac}
R.~B. Dingle.
\newblock On the {Fermi-Dirac} integrals {${\cal F}_p(\eta) = (p!)^{-1}
  \int_0^\infty \varepsilon^p (e^{\varepsilon-\eta}+1)^{-1} d\varepsilon$}.
\newblock {\em Applied Scientific Research}, 6(4):225--239, 1957.

\bibitem{Fisher}
Michael~E. Fisher.
\newblock Walks, walls, wetting and melting.
\newblock {\em Journal of Statistical Physics}, 34(5--6):667--729, 1984.

\bibitem{galluccio-loebl:pfaffian}
Anna Galluccio and Martin Loebl.
\newblock On the theory of {Pfaffian} orientations. {I}.\ {Perfect} matchings
  and permanents.
\newblock {\em Electronic Journal of Combinatorics}, 6(1):\#R6, 1999.

\bibitem{huang-wu-kunz:5v}
H.~Y. Huang, F.~Y. Wu, H.~Kunz, and D.~Kim.
\newblock Interacting dimers on the honeycomb lattice: {An} exact solution of
  the five-vertex model.
\newblock {\em Physica \textrm{A}}, 228(1--4):1--32, 1996.
\newblock arXiv:cond-mat/9510161.

\bibitem{kasteleyn:pfaffian}
P.~W. Kasteleyn.
\newblock Graph theory and crystal physics.
\newblock In Frank Harary, editor, {\em Graph Theory and Theoretical Physics}.
  Academic Press, 1967.

\bibitem{kenyon:dimers}
Richard Kenyon.
\newblock Local statistics of lattice dimers.
\newblock {\em Annales de l'Institut Henri Poincar\'e -- Probabilit\'es et
  Statistiques}, 33(5):591--618, 1997.
\newblock arXiv:math.CO/0105054.

\bibitem{lewin:polylogs}
Leonard Lewin.
\newblock {\em Polylogarithms and Associated Functions}.
\newblock North Holland, 1981.

\bibitem{lindelof:residus}
Ernst Lindel{\"o}f.
\newblock {\em Le Calcul des R\'esidus et ses Applications a la Th\'eorie des
  Fonctions}.
\newblock Gauthier-Villars, 1905.

\bibitem{lu-wu:dimer1}
W.~T. Lu and F.~Y. Wu.
\newblock Dimer statistics on the {M{\"o}bius} strip and the {Klein} bottle.
\newblock {\em Physics Letters \textrm{A}}, 259(2):108--114, 1999.
\newblock arXiv:cond-mat/9906154.

\bibitem{lu-wu:dimer2}
Wentau~T. Lu and F.~Y. Wu.
\newblock Close-packed dimers on nonorientable surfaces, 2001.
\newblock arXiv:cond-mat/0110035.

\bibitem{Nagle}
J.~F. Nagle.
\newblock Theory of biomembrane phase transitions.
\newblock {\em Journal of Chemical Physics}, 58:252--264, 1973.

\bibitem{pickard:polylogs}
W.~F. Pickard.
\newblock On polylogarithms.
\newblock {\em Publicationes Mathematicae}, 15(1--4):33--43, 1968.

\bibitem{popkov-kim-huang-wu:3d}
V.~Popkov, Doochul Kim, H.~Y. Huang, and F.~Y. Wu.
\newblock Lattice statistics in three dimensions: solution of layered dimer and
  layered domain wall models.
\newblock {\em Physical Review \textrm{E}}, 56(4):3999--4008, 1997.
\newblock arXiv:cond-mat/9703065.

\bibitem{Spohn}
M.~{Pr\"ahofer} and H.~Spohn.
\newblock An exactly solved model of three-dimensional surface growth in the
  anisotropic {KPZ} regime.
\newblock {\em Journal of Statistical Physics}, 88(5--6):999--1012, 1997.
\newblock arXiv:cond-mat/9612209.

\bibitem{regge-zecchina:pfaffian}
Tullio Regge and Riccardo Zecchina.
\newblock Combinatorial and topological approach to the {3D} {Ising} model.
\newblock {\em Journal of Physics \textrm{A}}, 33:741--761, 2000.
\newblock arXiv:cond-mat/9909168.

\bibitem{Sheffield}
Scott Sheffield.
\newblock Gibbs measure uniqueness results for {Lipschitz} random surfaces,
  2001.
\newblock In preparation.

\bibitem{Tabachnikov}
Serge Tabachnikov.
\newblock {\em Billiards}.
\newblock Panoramas et Synth\`eses. Soci\'et\'e Math\'ematique de France, 1995.

\bibitem{tesler:pfaffian}
Glenn Tesler.
\newblock Matchings in graphs on non-orientable surfaces.
\newblock {\em Journal of Combinatorial Theory, Series B}, 78:198--231, 2000.

\bibitem{truesdell:polylogs}
C.~Truesdell.
\newblock On a function which occurs in the theory of the structure of
  polymers.
\newblock {\em The Annals of Mathematics}, 46(1):144--157, 1945.

\bibitem{Wu67}
F.~Y. Wu.
\newblock Exactly soluable model of the ferroelectric phase transition in two
  dimensions.
\newblock {\em Physical Review Letters}, 18(15):605--607, 1967.

\bibitem{Wu68}
F.~Y. Wu.
\newblock Remarks on the modified potassium dihydrogen phosphate model of a
  ferroelectric.
\newblock {\em Physical Review}, 168(2):539--543, 1968.

\bibitem{WH:93}
F.~Y. Wu and H.~Y. Huang.
\newblock Exact solution of a vertex model in $d$ dimensions.
\newblock {\em Letters in Mathematical Physics}, 29(3):205--213, 1993.

\bibitem{WH}
F.~Y. Wu and H.~Y. Huang.
\newblock Exact solution of a lattice model of flux lines in superconductors.
\newblock {\em Physika \textrm{A}}, 205:31--40, 1994.

\end{thebibliography}
\bibliographystyle{plain}

\enlargethispage*{3\baselineskip}

\appendix
\section{Polylogarithms}\label{sec:polylogs}

The polylogarithm function $\Li_\nu(z)$ is defined by
$$\Li_\nu(z)=\sum_{n=1}^{\infty}\frac{z^n}{n^\nu}$$ for $|z|<1$ and by
analytic continuation elsewhere.  We may for instance write
$\Li_1(z)=-\log(1-z)$, $\Li_0(z) = z/(1-z)$, and $\Li_{-1}(z) = z/(1-z)^2$.
The polylogarithm function has the convenient property
that $$z \frac{\partial \Li_\nu(z)}{\partial z} = \Li_{\nu-1}(z).$$
Appell's integral expression
\begin{equation}\label{Liint}
\Li_\nu(z) = \frac1{\Gamma(\nu)}\int_0^\infty\frac{z s^{\nu-1}\,ds}{e^s-z}
\end{equation}
is valid for $\Re(\nu)>0$ and $z\not\in[1,\infty)$, and defines the
principal branch of the polylogarithm.  The polylogarithm has an
interesting Riemann surface.  When one crosses the branch cut
$[1,\infty)$ in the positive direction, the polylogarithm increases by
$\frac{2\pi i}{\Gamma(\nu)} (\log z)^{\nu-1}$.  (Both this and the
defining series expansion for $\Li_\nu$ are readily derived from Appell's
integral expression.)  For nonpositive integer $\nu$ this quantity is
$0$, consistent with the fact that $\Li_\nu$ is a rational function
for these values of $\nu$.  Unless $\nu$ is an integer $\leq1$, the
$\frac{2\pi i}{\Gamma(\nu)} (\log z)^{\nu-1}$ term creates a second
branch point at $z=0$ off the principal branch.  The function
$\Li_\nu(z)$ is analytic in both $z$ and $\nu$, except for a
singularity at $z=1$ (and $z=0$ off the principal branch).
For further background see Bateman and Erd\'elyi \textit{et al.\/}
\cite[Chapt.~1 \S11]{bateman-erdelyi:vol1},
Truesdell \cite{truesdell:polylogs},
Dingle \cite{dingle:fermi-dirac,dingle:bose-einstein},
and Lewin \cite{lewin:polylogs}.

When $z$ is on the principal branch and near $1$, the series expansion
\begin{equation}
\label{Li-pole}
\Li_\nu(z) = \Gamma(1-\nu) (-\log z)^{\nu-1} + \sum_{n=0}^\infty \zeta(\nu-n)\frac{(\log z)^n}{n!}
\end{equation}
was given by Lindel\"of \cite[pp.\ 138--141]{lindelof:residus}
(derivations are also given in \cite{truesdell:polylogs}
and \cite{chang-peres:zeta}),
and is absolutely convergent when $|\log z|<2\pi$.

We also use the asymptotic series expansions for $\Li_\nu(z)$ when $|z|\rightarrow\infty$.  For large positive $x$ these are
\begin{equation}
\label{Li(-x)}
\Li_\nu(-x) = -\cos(\pi\nu) \Li_\nu(-1/x) +
 2\sum_{k=0}^{\infty}
  \frac{\left(1/2^{2k-1}-1\right) \zeta(2 k)}{\Gamma(\nu+1-2k)} (\log x)^{\nu-2k}
\end{equation}
and
\begin{equation}
\label{Li(x)}
\Li_\nu(x) = -\cos(\pi\nu) \Li_\nu(1/x)
\pm \pi i \frac{(\log x)^{\nu-1}}{\Gamma(\nu)}
+
 2\sum_{k=0}^{\infty}
  \frac{\zeta(2 k)}{\Gamma(\nu+1-2k)} (\log x)^{\nu-2k}.
\end{equation}
Aside from the $-\cos(\pi\nu)\Li_\nu(\pm1/x)$ terms, the asymptotic
series expansions \eqref{Li(-x)} and \eqref{Li(x)} were derived by
Sommerfeld and Clunie respectively.  Rhodes first derived \eqref{Li(-x)}
for the case of integer $\nu$, and for these $\nu$ the expansions
\eqref{Li(-x)} and \eqref{Li(x)} have only finitely many nonzero terms.
For noninteger $\nu$ the expansions \eqref{Li(-x)} and \eqref{Li(x)} diverge,
and the $-\cos(\pi\nu)\Li_\nu(\pm1/x)$ term is dominated by each term in the
divergent series expansion, so one may wonder what role it plays.
Dingle \cite[eqn~17]{dingle:fermi-dirac} \cite[eqn~11]{dingle:bose-einstein}
showed how to make practical computational
use of these series by truncating the series after finitely many terms
and providing convergent series expansions for the remainder when the
$-\cos(\pi\nu)\Li_\nu(\pm1/x)$ term is present.  For positive $x>1$ we
are evaluating $\Li_\nu(x)$ on the branch cut of the principal branch;
the $\pm$ sign is positive when the branch cut is just above the real
axis, which is the usual convention.
See also Pickard \cite[eqn 3.5]{pickard:polylogs} for the asymptotics
for large complex values of $z$.




\section{Notation}\label{sec:notation}

$a$, $b$, $c$: weights of edges in the three different directions.  $a$ is weight for vertex not being in a loop, $b$ is weight for horizontal step, and $c$ is weight for vertical step.

\noindent
$m$: horizontal length of torus

\noindent
$n$: vertical length of torus

\noindent
Simplifying assumptions: $b<a$, $b/a=\Theta(1)$, $1-b/a=\Theta(1)$, $n=\Theta(m)$

\noindent
$p/q$: rational approximation of $nb/(mc)$, $\gcd(p,q)=1$.
Intuitively $p$ is the number of horizontal windings of loops,
and $q$ is the number of vertical windings of loops, but this
intuition is not quite accurate if ``ratcheting'' takes place.

\noindent
$W$ is a measure of close the rational approximation $p/q \approx nb/(mc)$
needs to be.  For fixed $p$ and $q$, $W=\Theta(1/\sqrt n)$.  More precisely,

\noindent
$W = \sqrt{q \eps}/(\pi p) = \sqrt{2 q n a b}/(p m (a-b)) \approx \sqrt{2/(pmc)} \approx \sqrt{2/(qnb)}$

\noindent
$\alpha$ is measure of the error in the approximation $p/q \approx nb/(mc)$ on the scale of $W$:

\noindent
$\displaystyle \frac{n b q}{m c p} \approx \frac{n b q}{m (a-b) p} = 1+\alpha W.$

\noindent
$A = \left(\frac{c}{a-b}\right)^n$

\noindent
$A$ is a measure of how close the weights $a$, $b$, and $c$ are to the phase transition $a=b+c$.  When $a=b+c$, $A=1$, and $A$ is sensitive to perturbations on the order of $\Theta(1/n)$.  $A^q$ appears in many of our formulas, and is approximately symmetric in the parameters: $$ \log A^q \approx (b+c-a) \sqrt{\frac{m n p q}{b c}}$$

\noindent
$\phi = \frac{2\pi n b}{m(a-b)}$

\noindent
$\eps = \frac{2\pi^2 n a b}{m^2(a-b)^2}$

\noindent
$z_k=-e^{i\theta_k}$, $\theta_k = 2\pi k/m$; $k\in\Z_m+\sigma/2$ for $\Zst$

\noindent
$r_k=(a-b)^{n}|a+bz_k|^{-n} = e^{-\eps k^2+O(k^4/m^3)}$

\noindent
$\phi_k=\arg(a+bz_k)^{-n} = \phi k+O(k^3/m^2)$

\end{document}